\magnification=\magstep1
\input amssym.def

\baselineskip=1.2\baselineskip

\def\card{{{\mathop{\rm{Card}}}}}
\def\supp{{\rm supp}\,}

\def\eps{\varepsilon}
\def\noi{\noindent}
\def\bone{{\bf 1}}
\def\R{{\Bbb R}}
\def\C{{\bf C}}
\def\N{{\Bbb N}}

\def\D{{\bf D}}

\def\cE{{\cal E}}
\def\cL{{\cal L}}
\def\calw{{{\cal W}}}

\def\wt{\widetilde}
\def\wh{\widehat}

\def\sqr#1#2{{\vcenter{\vbox{\hrule height.#2pt
        \hbox{\vrule width.#2pt height#1pt \kern#1pt
           \vrule width.#2pt}
        \hrule height.#2pt}}}}
\def\square{\mathchoice\sqr56\sqr56\sqr{2.1}3\sqr{1.5}3}
\def\build#1_#2^#3{\mathrel{\mathop{\kern 0pt#1}\limits_{#2}^{#3}}}
\def\un{\underline}
\def\ov{\overline}

\input epsf.tex
\newdimen\epsfxsize

\centerline{\bf SUPER-BROWNIAN MOTION }
\smallskip
\centerline{\bf WITH REFLECTING HISTORICAL PATHS}
\footnote{$\empty$}{\rm Research partially supported by NSF grant DMS-9700721.}
\bigskip
\centerline{\bf Krzysztof Burdzy}
\centerline{\bf Jean-Fran\c cois Le Gall}
\vskip1truein

{\bf Abstract}. We consider super-Brownian motion
whose historical paths reflect from each other, unlike those
of the usual historical super-Brownian motion.
We prove tightness for the family of
distributions corresponding to a sequence of discrete approximations
but we leave the problem of uniqueness of the limit open.
We prove a few results about path behavior for processes
under any limit distribution. In particular, we show
that for any $\gamma>0$,
a ``typical'' increment of a reflecting historical path
over a small time interval $\Delta t$ is not greater
than $(\Delta t)^{3/4 - \gamma}$.

\bigskip

\noindent{\bf 1. Introduction}.
\medskip
The present article has been inspired
by two probabilistic models---superprocesses
with interactions and reflected particle systems.

The first person to study a reflecting system
of particles was
Harris [H] who considered an infinite system of Brownian
particles on the line.
He proved that if the initial positions of the particles
are points of a Poisson point process, then for a large
time $t$ the distribution of a single particle is normal
with the standard deviation $(2t/\pi)^{1/4}$. Spitzer [S] analyzed
a similar model with particles moving along
straight lines between collisions. See
[DGL1, DGL2, G, Ho] for related results.

The simplest superprocesses, for example, super-Brownian
motion, are continuum limits of branching systems
in which the branching mechanism is independent
of the positions of particles. There has been
considerable activity studying models with interactions.
Many articles are devoted to models with catalysts,
see, e.g., [DF, De]. Various other models
with interactions are discussed in [AT, BHM, EP, P3].
See in particular [P4] and references therein.

We will study a model similar to that introduced by Harris,
in that we will start with linear Brownian motion as the spatial
process.
We will attempt to build a corresponding superprocess
with historical paths that do not cross over although they may touch
each other.

Our construction is based on a sequence of discrete
approximations. Consider for every $\varepsilon\in(0,1]$
a branching particle system which starts initially
with $N_\varepsilon$ particles located respectively at
$x^\varepsilon_1\leq\cdots\leq x^\varepsilon_{N_\varepsilon}$.
Particles move independently in space according to linear
Brownian motion and are subject to
critical binary branching at rate $\varepsilon^{-1}$.
To be specific, the lifetimes of the particles are exponential
with parameter $\varepsilon^{-1}$ and
when a particle dies it gives rise to $0$ or $2$ new particles
with probability $1/2$.

Let us now introduce our basic assumptions. Let
$$\mu_\varepsilon:=\varepsilon\sum_{j=1}^{N_\varepsilon}
\delta_{x^\varepsilon_j}
$$
and assume that there is a finite measure $\mu$ on $\R$ such that
$$\mu_\varepsilon
\build{\hbox to 8mm{\rightarrowfill}}_{\varepsilon\to 0}^{\rm (w)} \mu,
\eqno{(1.1)}$$
where the notation (w) indicates weak convergence in the space
$M_f(\R)$ of finite measures on $\R$. In addition, if
$\supp\mu$ denotes the topological support of $\mu$, we assume that
$$\supp\mu_\varepsilon
\build{\hbox to 8mm{\rightarrowfill}}_{\varepsilon\to 0}^{} \supp\mu,
\eqno{(1.2)}$$
in the sense of the Hausdorff metric
on compact subsets of $\R$ (in particular, we assume that $\supp\mu$
is compact).

 Let $X^\varepsilon_t$ denote the
random measure equal to $\varepsilon$ times the sum of the
Dirac point masses at the positions of particles alive at time $t$.
Then,
$$(X^\varepsilon_t,t\geq 0)
\build{\hbox to 8mm{\rightarrowfill}}_{\varepsilon\to 0}^{\rm
(d)} (X_t,t\geq 0),\eqno{(1.3)}$$
where the limit process is super-Brownian motion in $\R$ with branching rate
$\gamma=1$ (throughout this work we consider only this branching rate)
and initial value $\mu$,
and the convergence holds in distribution in the
Skorohod space $\D(\R_+,M_f(\R))$. The convergence (1.3) is  the
standard approximation of super-Brownian motion (see
e.g. [P4]). Note that assumption (1.2) is not
needed for (1.3) but it guarantees that the graph of $X^\varepsilon$
also converges in distribution to the graph of $X$ (see Lemma 2.3 below),
a property that plays an important role in our arguments.

For each particle alive at time $t$, we can consider
its historical path, which is the element of $\C([0,t],\R)$
obtained by concatenating the trajectories of the ancestors
of the given particle up to time $t$. Denote by $Y^\varepsilon_t$
the historical measure equal to $\varepsilon$ times
the sum of the Dirac point masses at
the historical paths of the particles alive at time $t$ ($Y^\varepsilon_t$
is thus a random measure on the set $\C([0,t],\R)$
of continuous mappings from $[0,t]$ into $\R$).
Then the convergence (1.3) can be reinforced as
$$(Y^\varepsilon_t,t\geq 0)
\build{\hbox to 8mm{\rightarrowfill}}_{\varepsilon\to 0}^{\rm
(d)} (Y_t,t\geq 0),\eqno{(1.4)}$$
where the limit process is now historical super-Brownian motion started
at $\mu$.

For every $\varepsilon>0$, we can use the original
branching particle system to construct a new system with reflection.
The branching mechanism (critical binary branching at rate
$\varepsilon^{-1}$) is the same as in the original system, but the particle
paths
in the new system reflect against each other.
A precise construction is given in Section 3, but let us
give an informal description. The reflected system
is such that for every $t\geq 0$, the set of
positions of particles at time $t$ is the same as in
the original system, and in particular the branching times
are the same. During the time interval between $0$ and the first
branching time, the vector of positions of the particles
labeled $1,2,\ldots,N_\varepsilon$ in the reflected system
is the increasing rearrangement of the vector of positions
of the particles in the original system. Suppose that at the
first branching time, denoted by $\xi$, a particle dies and gives rise to
2 children. If the location of this particle is the $j$-th coordinate
in the increasing rearrangement of the vector of positions at time
$\xi-$, we will say that in the reflected system particle $j$ has
given rise to two children labeled $j1$ and $j2$. Then on the
interval between $\xi$ and the second branching time, the
vector of positions of the particles labeled
$1,\ldots,j-1,j1,j2,j+1,\ldots,N_\varepsilon$
in the reflected system is again the increasing rearrangement of the vector
of positions
of the particles in the original system. We can easily continue this
construction
by induction.

Denote
by $\wt X^\varepsilon_t$ and $\wt Y^\varepsilon_t$ the analogues
of $X^\varepsilon_t$ and $Y^\varepsilon_t$ for the
the system with reflection. We have
$\wt X^\varepsilon_t=X^\varepsilon_t$ since the set of positions of
particles is the same at every
time
$t$ in both systems. On the other hand, $\wt Y^\varepsilon_t$ is typically
very different
from $Y^\varepsilon_t$. Indeed, the following  property holds for any two
paths $w$, $w'$ in the
support of $\wt Y^\varepsilon_t$: Either $w(r)\leq w'(r)$ for every $0\leq
r\leq t$, or
$w(r)\geq w'(r)$ for every $0\leq r\leq t$.

The main purpose of this work is to try to understand the
limiting behavior of the branching particle system with
reflection as $\varepsilon\to 0$. Our primary objective was to get
an analogue of the convergence (1.4) when the processes
$Y^\varepsilon$ are replaced by $\wt Y^\varepsilon$,
giving information about the individual paths in the
system with reflection. We did not completely succeed in this task, but
we can prove the following result, where ${\cal W}$ denotes the set
of all stopped paths, or equivalently the union over all $t\geq 0$
of the sets $\C([0,t],\R)$.
\vfill\eject
\noi{\bf Theorem 1.1}. {\it Let ${\cal E}$ be a sequence of positive numbers
converging to $0$. The laws of the
processes $(\wt Y^\varepsilon_t,t\geq 0)$ for $\varepsilon\in{\cal E}$
are tight in the space of all probability measures on
the Skorohod space $\D([0,\infty),M_f({\cal W}))$. Furthermore any
limiting distribution is supported on $\C([0,\infty),M_f({\cal W}))$. }
\medskip
Hence, by extracting a subsequence if necessary, we can assume that
the sequence of processes $\wt Y^\varepsilon$ converges in
distribution towards a process $\wt Y$ with continuous paths with
values in $M_f({\cal W})$. Note that, for every $t\geq 0$, the measure
$\wt Y_t$ is supported on $\C([0,t],\R)$. Although the question of
uniqueness of the limit
remains unsolved,  we
are able to derive
several results on the path behavior of the process $\wt Y$.

First note that, since $\wt X^\varepsilon_t=X^\varepsilon_t$
for every $t\geq 0$, the convergence (1.1) implies that the
$M_f(\R)$-valued process
$\wt X$ defined by
$$\langle\wt X_t,\varphi\rangle=\int \wt Y_t(dw)\,\varphi(w(t))$$
is a super-Brownian motion started at $\mu$. In particular, it is known
(see [KS], [R]) that a.s. for every $t>0$ the measure $\wt X_t(dy)$ has a
density
denoted by $x_t(y)$, and that there exists
a jointly continuous modification of $(x_t(y),t>0,y\in \R)$.

The next result shows that for any $\gamma>0$,
a typical oscillation of a reflecting
historical path is not greater than $(\Delta t)^{{3\over 4}-\gamma}$,
and hence much smaller than a typical Brownian oscillation
$(\Delta t)^{{1\over 2}}$. This result is consistent with
the Harris [H] estimate, if we translate the large-time
asymptotics to small-time asymptotics.
\medskip
\noi{\bf Theorem 1.2}. {\it Almost surely for every $t>0$ and every
$r\in(0,t)$,
for every path $w\in\supp\wt Y_t$, the condition $x_r(w(r))>0$ implies
that, for every $\gamma>0$,
$$\limsup_{\delta\downarrow 0} {|w(r+\delta)-w(r)|\over \delta^{{3\over
4}-\gamma}}=0.$$}

A more precise version of Theorem 1.2 is given in Section 5 (Theorem 5.10).
It is not hard to check that if we fix
$t>0$ and $r\in (0,t)$ (fixing $r$ is in fact enough), the
condition $x_r(w(r))>0$, and thus the conclusion of the theorem, will
hold for every path $w\in\supp\wt Y_t$, a.s. Alternatively, for every fixed
$t>0$, the
conclusion of Theorem 1.2 holds for a set of values of $r\in(0,t)$
of full Lebesgue measure, for every $w\in\supp\wt Y_t$.
We believe that
$\delta^{3\over 4}$ is the ``typical'' size for the
oscillation of a historical reflected path
although we have no lower bound justifying this claim.

We also study the behavior of reflected historical paths at a branching
point. If $w$ and $w'$ are two reflected historical paths that
coincide up to time $r>0$ (meaning informally that the
corresponding ``particles'' have the same ancestor up to
time $r$), we show that the distance between $w(r+\delta)$
and $w'(r+\delta)$ grows linearly as a function of $\delta$, up
to logarithmic corrections. The precise statement is as follows.
\medskip
\noi{\bf Theorem 1.3}. {\it Let $t>0$. If $w$ and $w'$ are
two distinct elements of $\C([0,t],\R)$, we set
$$\gamma_{w,w'}=\inf\{r\geq 0:w(r)\not =w'(r)\}.$$
Then a.s. for any two distinct paths $w,w'\in\supp\wt Y_t$
such that $\gamma_{w,w'}>0$, we have
$$\limsup_{\delta\downarrow 0}
{|w(\gamma_{w,w'}+\delta)-w'(\gamma_{w,w'}+\delta)|\over
2\delta\log|\log\delta
|}=x_{\gamma_{w,w'}}(w(\gamma_{w,w'}))>0$$
and, for every $\gamma>0$,
$$\lim_{\delta\downarrow 0} {|w(\gamma_{w,w'}+\delta)-w'(\gamma_{w,w'}+\delta)|
\over \delta|\log\delta|^{-1-\gamma}}=\infty.$$}

Our proofs rely on several known results on super-Brownian
motion.
In particular, we use the Brownian snake idea [L2]
in an essential way, both in the proofs and
for giving more precise versions of the results.
For instance, as a key step towards Theorem 1.1, we
get a uniform continuity result (Theorem 4.1) for the
historical paths of the approximating branching particle systems
with reflection. The proof of this result requires some
precise information about the genealogical structure
of the approximating systems, which seems to be more
easily accessible via the snake approach (cf Lemma 2.1 below).

For an introduction to the theory of superprocesses
(measure-valued diffusions) and historical processes,
the reader may consult [Da, Dy, DP, L2, P4].

The paper is organized as follows. Section 2 describes the
specific coding that we use to represent the genealogical
structure of the approximating branching particle systems.
This section also contains a few important preliminary
results. Section 3 presents the construction of
the systems with reflection. Tightness results are given
in Section 4, including a more precise form of Theorem 1.1.
Section 5 contains the proof of Theorem 1.2, and is the
most technical part of the paper. Finally, Theorem 1.3 is
proved in Section 6.

We are grateful to Carl Mueller, Ed Perkins, Tokuzo Shiga
and Roger Tribe for very useful advice.

\bigskip

\noindent{\bf 2. Coding discrete trees}.
\medskip
We will describe a method that provides a coding of
the genealogy of the branching particle systems introduced in Section 1,
in a consistent way for all values of the parameter $\varepsilon\in(0,1]$.
This method involves embedding branching trees in a path of reflected Brownian
motion, and is based on [L1] (see also [NP]).
\medskip
\noindent{\bf 2.1 \ Markov chains embedded in reflected Brownian motion}.
\medskip
Let $\beta=(\beta_s, s\ge 0)$ be distributed as twice a reflected Brownian
motion on
$\R_+$:
$$(\beta_s,s\geq 0)\buildrel{\rm(d)}\over = (2|B_s|,s\geq 0),$$
where $B$ is a standard linear Brownian motion, with $B_0=0$. The reason for
the factor $2$ will be clear later. We denote by $(L^x_s,x\geq 0,s\geq 0)$ the
jointly continuous family of local times of $\beta$, normalized in such
a way that,
for every nonnegative Borel function $\varphi$ on $\R_+$,
$$\int_0^t \varphi(\beta_s)\,ds=\int_{\R_+} \varphi(x)L^x_tdx.$$
Also set $\tau_r=\inf\{s\geq 0:L^0_s>r\}$, for every $r>0$.

For every $\eps\in(0,1]$, we introduce a sequence of
stopping times $(T^\eps_k, k=0, 1, \ldots)$ defined inductively as
follows:
$$\eqalign{T^\eps_0&=\inf\{s\ge 0: \beta_s=2\eps\},\cr
T^\eps_{2k+1}&=\inf\{u\ge T^\eps_{2k}: \sup_{T^\eps_{2k}\le s\le
u}\beta_s-\beta_u=2\eps\},\cr
T^\eps_{2k+2}&=\inf\{u\ge T^\eps_{2k+1}: \beta_u-
\mathop{\inf}\limits_{T^\eps_{2k+1}\le s\le u}\beta_s=2\eps\}.\cr}$$
It is simple to check that the
variables $T^\eps_0, T^\eps_1-T^\eps_0, T^\eps_2-T^\eps_1,
\ldots$ are independent and identically distributed. To see
this, note that if $(\gamma_t, t\ge 0)$ is a reflected Brownian motion with
initial value $\gamma_0=b\ge 0$, the process
$$\gamma_t-\mathop{\inf}\limits_{0\le s\le
t}\gamma_s$$
is again a reflected Brownian motion, with initial value $0$, and also
observe that $\beta_{T^\eps_{2k}}\ge 2\eps$ for every $k$.

As $E(T^\eps_0)=\eps^2$, standard arguments  show that for every
$K>0$
$$\sup_{s\le K}\left|T^\eps_{[s/\eps^2]}-s\right|
\mathop{\longrightarrow}\limits^{\mathop{\rm{a.s.}}}_{\eps\to 0}0.
\eqno{(2.1)}$$
(First establish this convergence along the sequence $\varepsilon_n=n^{-2}$
and then use monotonicity arguments.) Thus,
$$\sup_{s\le K}\left|\beta_{T^\eps_{[s/\eps^2]}}-\beta_s\right|
\mathop{\longrightarrow}\limits^{\mathop{\rm{a.s.}}}_{\eps\to 0}0.
\eqno{(2.2)}$$

For $k=0, 1, \ldots$, set
$$\eqalign{S^\eps_{2k}&=\beta_{T^\eps_{2k}}-2\eps,\cr
S^\eps_{2k+1}&=\beta_{T^\eps_{2k+1}}.\cr}$$
It is easy to verify that $(S^\eps_k, k=0, 1, 2, \ldots)$ is a
time-inhomogeneous
Markov chain with values in $\R_+$, whose law can be described as follows
(see [L1] Section 3 for details):
$S^\eps_0=0$ and $S^\eps_{2k+1}$ has the same distribution as
$S^\eps_{2k}+U$, where $U$ is an exponential variable with mean
$2\eps$, independent of $S^\eps_{2k}$, $S^\eps_{2k+2}$ has the
same distribution as $(S^\eps_{2k+1}-V)_+$ where $V$ is exponential with
mean $2\eps$, independent of $S^\eps_{2k+1}$.

{}From (2.2), we have a.s.\ for every $K>0$,
$$\sup_{s\le K}\left|S^\eps_{[s/\eps^2]}-\beta_s\right|
\mathop{\longrightarrow}\limits^{\mathop{\rm{a.s.}}}_{\eps\to 0}0.$$
We then define a continuous-time process $(\beta^\eps_s, s\ge 0)$ by
setting
$$\beta^\eps_{k\eps^2}=S^\eps_k\hbox{ for }k=0, 1, 2, \ldots$$
and by interpolating linearly on intervals of the form $[k\eps^2,
(k+1)\eps^2]$. It is obvious that we also have
$$\sup_{s\le K}\left|\beta^\eps_s-\beta_s\right|
\mathop{\longrightarrow}\limits^{\mathop{\rm{a.s.}}}_{\eps\to
0}0.\eqno{(2.3)}$$

\noi{\bf 2.2 \ The correspondence between excursions and trees}
\medskip
With each excursion of $\beta^\varepsilon$ away from $0$, we can
associate a marked tree representing the genealogical structure of a
Galton-Watson branching process with critical binary branching
at rate $\varepsilon^{-1}$, starting with one individual
(the ancestor) at time $0$. Here a marked tree consists
of the set $\cal T$ of edges (i.e., particles), which is a subset of
$${\bf U}:=\bigcup_{n=0}^\infty \{1,2\}^n
\qquad\hbox{(by convention, }\{1,2\}^0=\{\emptyset\}),$$
and the family
$(\ell_u,u\in{\cal T})$ of lengths of edges (i.e., lifetimes of
particles).

\bigskip
\vbox{
\epsfxsize=5.0truein
  \centerline{\epsffile{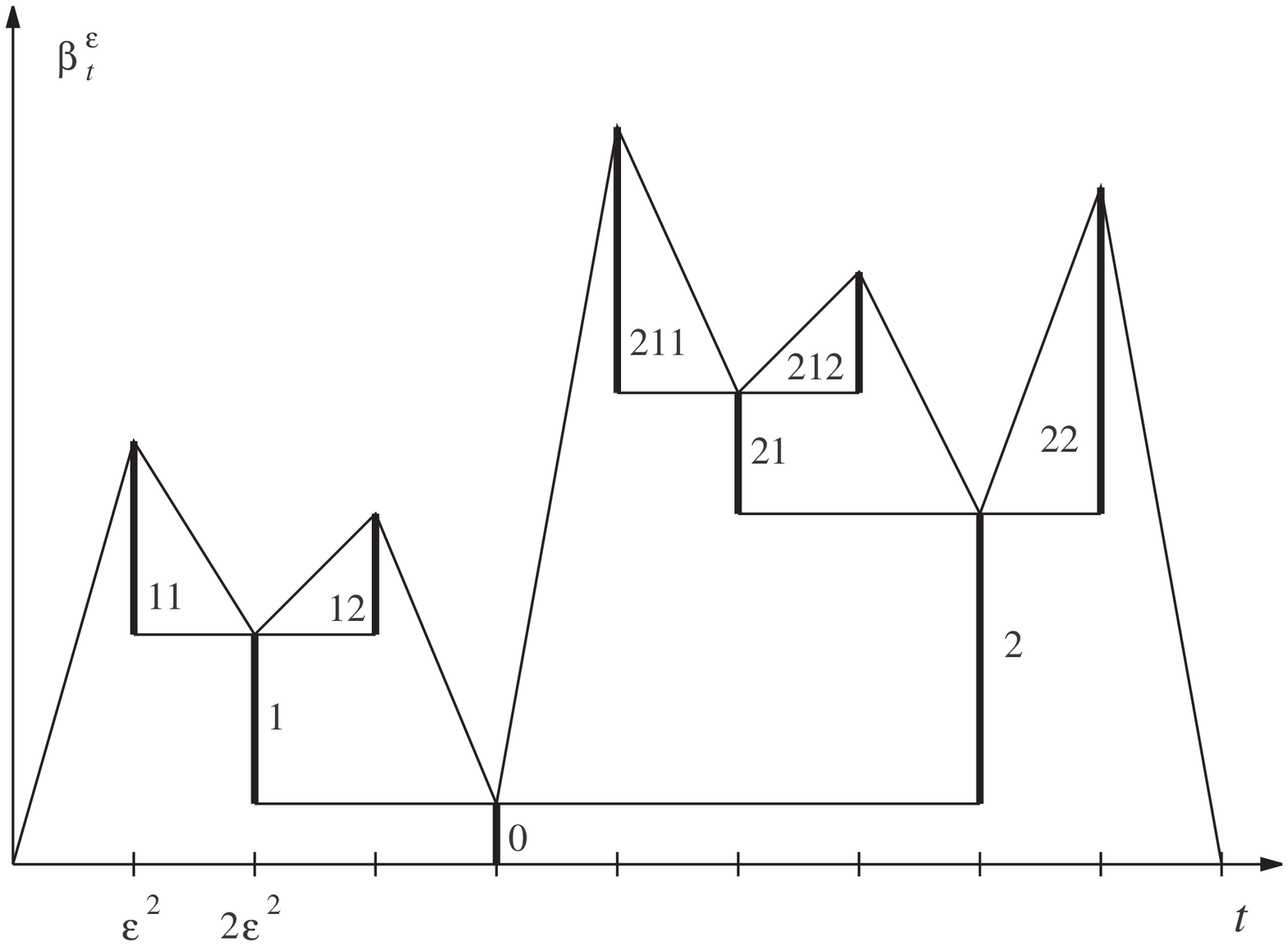}}

\centerline{Figure 1.}
}
\bigskip

This
correspondence is explained in
Fig.~1 for the first excursion of $\beta^\varepsilon$
away from $0$. Informally, if
$(i\varepsilon^2,j\varepsilon^2)$ is the interval corresponding to an
excursion of $\beta^\varepsilon$, the
lifetime $\ell_\emptyset$ of the individual at the root
of the associated tree is the minimum of
$\beta^\varepsilon$ over
$[(i+1)\varepsilon^2,(j-1)\varepsilon^2]$ and this individual has two children
if and only if
$j-1>i+1$. In that case, by decomposing the excursion
restricted to $[(i+1)\varepsilon^2,(j-1)\varepsilon^2]$ at the time of its
minimum over this interval, we get two new excursions, each  of which codes
the
genealogical structure of
descendants of one of the ancestor's children. The construction of the
tree is then completed by induction. Note that each time
of the form $k\varepsilon^{-2}$ in the interval
$(i\varepsilon^2,j\varepsilon^2)$
corresponds to one edge of the tree (for instance the time of the minimum
over $[(i+1)\varepsilon^2,(j-1)\varepsilon^2]$ corresponds to $\emptyset$,
see Fig.~1). We refer to [L1]
Section 2 for a more precise description and a proof that this construction
yields the
family tree of a
Galton-Watson branching process with critical binary branching at rate
$\varepsilon^{-1}$.
(We can now explain  the factor $2$ in the definition
of $\beta$: We want the branching rate to be
$\varepsilon^{-1}$ and not $(\varepsilon/2)^{-1}$.)

There is a one-to-one correspondence between excursions of
$\beta^\varepsilon$ away from
$0$ and excursions of $\beta$ away from $0$ with height greater than
$2\varepsilon$:
If $k\varepsilon^2$ is the beginning of an excursion of $\beta^\varepsilon$,
then $T^\varepsilon_{k}$ is the hitting time of $2\varepsilon$
by the corresponding excursion of $\beta$. As in Section 1, consider
for every $\varepsilon\in(0,1]$ an integer $N_\varepsilon\geq 1$
and assume that the family $(\varepsilon
N_\varepsilon,\varepsilon\in(0,1])$
is bounded and that $\varepsilon N_\varepsilon$
converges to $a\geq 0$ as $\varepsilon\to 0$ (this follows
from (1.1) with $a=\langle \mu,1\rangle$). Let $\tau^\varepsilon$ denote
the $N_\varepsilon$-th return of $\beta^\eps$ to $0$.
{}From the previous observations, (2.1) and the standard
approximation of Brownian local times by upcrossing numbers, we have
$$\lim_{\eps\to 0}\tau^\eps=\tau_a,\hbox{ a.s.}$$
We will write $\tau=\tau_a$ for simplicity.

On the time interval $[0, \tau^\eps]$, the process $\beta^\eps$ makes
$N_\varepsilon$ independent excursions away from $0$. These excursions
can be viewed as representing the genealogical structure of the branching
particle system introduced in Section 1.
The set of edges, denoted by ${\cal T}_\varepsilon$, is then a random subset of
$\{1,\ldots,N_\varepsilon\}\times {\bf U}$ and
conditionally on ${\cal T}_\varepsilon$, the corresponding lengths are
independent exponentials with mean $\varepsilon$. The function
$(\beta^\varepsilon_s,s\in[0,\tau^\varepsilon])$ can be reconstructed from this
collection of marked trees as shown by Fig.~1. Notice that for this
reconstruction to be possible, it is
essential to order the trees and the different edges of every single tree.

\medskip

\noi {\bf 2.3 \ Discrete and continuous local times}
\medskip
One reason for considering the processes $\beta^\varepsilon$
comes from their relation with the upcrossing numbers of $\beta$.
We first define the (discrete) local times of $\beta^\varepsilon$. For every
$x\geq 0$ and $s\geq 0$, we define
$$L^{\varepsilon,x}_s=\varepsilon\,\card\{r\in[0,s):
\beta^\varepsilon_r=x\ {\rm and}\ \beta^\varepsilon_u>x\ {\rm
for}\ u\in(r,r+\delta],\ {\rm for\ some}\ \delta>0\}.$$
In other words, $\varepsilon^{-1}L^{\varepsilon,x}_s$ is the number of
upcrossings of $\beta^\varepsilon$ above level $x$ before time $s$.

Let $M^\varepsilon_s(x)$ denote the number of upcrossings
of $\beta$ from $x$ to $x+2\varepsilon$ completed before time $s$. More
precisely,
$M^\varepsilon_s(x)$ is the number of pairs $(u,v)$
with $0\leq u<v< s$, such that $\beta_u=x$, $\beta_r>x$ for
every $r\in(u,v)$ and $v=\inf\{r>u:\beta_r>x+2\varepsilon\}$.

Then, a.s. for every $x\geq 0$ and every integer $k\geq 1$, we have
$$L^{\varepsilon,x}_{(2k-1)\varepsilon^2}=L^{\varepsilon,x}_{2k\varepsilon^2}=
\varepsilon M^\varepsilon_{T^\varepsilon_{2k}}(x)
=\varepsilon M^\varepsilon_{T^\varepsilon_{2k-1}}(x).\eqno{(2.4)}$$
This identity is easily verified by induction
on $k$ (the sequence of stopping times
$(T^\varepsilon_k)$ was designed for this property to hold). See
also Proposition 7 of [L1].

\medskip
\noindent{\bf Lemma 2.1}. {\it We have a.s.
$$\lim_{\varepsilon\to 0}
\Big(\sup_{s\geq 0}\,\sup_{x\geq 0}|L^{\varepsilon,x}_{s\wedge\tau^\varepsilon}
-L^x_{s\wedge\tau}|\Big)=0.$$}
\smallskip
\noindent{\bf Proof}.
We first observe that a.s.
$$\lim_{\varepsilon\to 0}
\Big(\sup_{s\geq 0}\,\sup_{x\geq
0}|\varepsilon M^\varepsilon_{s\wedge\tau^\varepsilon} (x)-L^x_{s\wedge\tau}|
\Big)=0.\eqno{(2.5)}$$
For a fixed value of $x$, this is nothing but the classical approximation of
Brownian local time by upcrossing numbers, and excursion theory provides
precise estimates for the rate of convergence. Using these estimates and
monotonicity properties, it is then an easy task to
prove (2.5), i.e., the uniform version of the claim.

The statement of the lemma is now a simple
consequence of (2.1), (2.4) and (2.5).\hfill$\square$
\medskip
\noi{\bf Remark}. As an immediate consequence of Lemma 2.1 and the
joint continuity of Brownian local times, we get that
$$\lim_{\varepsilon,\delta\to 0}\Big(\;\sup_{s\geq 0}
\build{\sup_{x,x'\geq 0}}_{|x-x'|\leq \delta}^{}|L^{\varepsilon,x}_{s\wedge
\tau^\varepsilon}-L^{\varepsilon,x'}_{s\wedge \tau^\varepsilon}|
\Big)=0,\quad\hbox{a.s.}$$
Later, we will consider for every $\varepsilon\in(0,1]$ a process $\wt
\beta^\varepsilon$
with the same distribution as $\beta^\varepsilon$. If $\wt L^{\varepsilon,x}_s$
denote the discrete local times of $\wt \beta^\varepsilon$, the last
convergence
still holds in probability when $L^{\varepsilon,x}_s$
is replaced by $\wt L^{\varepsilon,x}_s$
(and $\tau^\varepsilon$ by $\wt\tau^\varepsilon$, with an
obvious notation).
\medskip
\noindent{\bf 2.4 \ Branching particle systems and discrete snakes}.
\medskip
We now consider the branching particle system of
Section 1, starting with $N_\varepsilon$ particles
located respectively at $x^\varepsilon_1, x^\varepsilon_2,
\ldots, x^\varepsilon_{N_\varepsilon}$. We may and will
assume that the genealogy of the descendants of particle $k$ (present at
$x^\varepsilon_k$ at time $0$) is given
by the tree associated with the $k$-th excursion of
$\beta^\varepsilon$ (cf subsection 2.2). We will refer to this
system as the $\varepsilon$-system of branching Brownian motions.

For our purposes, it will be convenient to view the
collection of paths
traced by the branching particles as the range of
a path-valued process called the discrete snake.

By definition, a stopped path is
a continuous mapping $w:[0,\zeta]\longrightarrow \R$,
where $\zeta=\zeta_w\geq 0$ is called the ``lifetime'' of $w$
(it is convenient to talk about the ``lifetime''
of a path although for technical reasons the path is stopped
rather than killed).
Let $\calw$ be the set of all stopped paths. Then
$\calw$
is a separable complete metric space for the distance
$$d(w,w')=|\zeta_w-\zeta_{w'}|+\sup_{t\geq 0}|w(t\wedge \zeta_w)-w'(t\wedge
\zeta_{w'})|.$$
For any $x\in\R$, we write $\un{x}$ for the trivial path
such that $\zeta_{\un{x}}=0$ and $\un{x}(0)=x$.
\smallskip
With every $s\in [0, \tau^\eps]$ we now associate a stopped path
$W^\eps_s\in\calw$
with lifetime $\beta^\varepsilon_s$. If $s\in [0, \tau^\eps)
\cap \varepsilon^2\N$ and $\beta^\eps_s=0$, then
$s$ is the starting time of the $k$-th excursion of
$\beta^\varepsilon$ away from $0$, for some
$k\in\{1,\ldots,N_\varepsilon\}$. We then
set $W^\eps_s=\underline{x}^\varepsilon_k$. For definiteness, we also set
$W^\varepsilon_{\tau^\varepsilon}=\underline{x}^\varepsilon_{N_\varepsilon}$. If
$s\in [0,
\tau^\eps)\cap\eps^2\N$ but $\beta^\eps_s>0$, we can associate with $s$ a
unique edge of the $k$-th tree, $k$ being the number of the excursion
straddling $s$. We then let $W^\eps_s$ be the historical
path of the particle in the system of branching Brownian motions that
corresponds to this edge. Notice that
the death time of this particle is  $\beta^\varepsilon_s$, and thus
$\zeta_{W^\varepsilon_s}=\beta^\varepsilon_s$. Finally if $s\in [0,
\tau^\eps]$ but $s\not\in\eps^2\N$,
we find an integer $j$ such that $j\eps^2<s<(j+1)\eps^2$, and let $l=j$ if
$\beta^\eps_{j\eps^2}>\beta^\eps_{(j+1)\eps^2}$, but $l=j+1$
if $\beta^\eps_{j\eps^2}\leq\beta^\eps_{(j+1)\eps^2}$. Then
we let $W^\eps_s$ be the path $W^\eps_{l\eps^2}$ stopped at time
$\beta^\eps_s$.

\smallskip
It is easy to see that conditionally on $(\beta^\varepsilon_s,s\geq 0)$
the process $(W^\varepsilon_{k\varepsilon^2},0\leq k\leq
\tau^\varepsilon/\varepsilon^2)$
is Markovian. To describe its conditional
distribution, let
$k\in\{0,\ldots,\tau^\varepsilon/\varepsilon^2\}$ and suppose that
$\beta^\varepsilon_{(k+1)\varepsilon^2}>0$ (otherwise
$W^\varepsilon_{(k+1)\varepsilon^2}
=\un{x}^\varepsilon_j$,
if $(k+1)\varepsilon^2$ is the
starting point of the $j$-th excursion of $\beta^\varepsilon$). If
$\beta^\varepsilon_{(k+1)\varepsilon^2}\leq
\beta^\varepsilon _{k\varepsilon^2}$ (which occurs if $k$ is odd) then
$W^\varepsilon_{(k+1)\varepsilon^2}$
is simply the restriction of $W^\varepsilon_{k\varepsilon^2}$ to $[0,
\beta^\varepsilon_{(k+1)\varepsilon^2}]$. On the other hand, if
$\beta^\varepsilon_{(k+1)\varepsilon^2} > \beta^\varepsilon
_{k\varepsilon^2}$, then  $W^\varepsilon_{(k+1)\varepsilon^2}$
is obtained from $W^\varepsilon_{k\varepsilon^2}$ by ``adding
at the tip of $W^\varepsilon_{k\varepsilon^2}$''
a Brownian path of length  $\beta^\varepsilon_{(k+1)\varepsilon^2} -
\beta^\varepsilon _{k\varepsilon^2}$ independent of
$(W^\varepsilon_{j\varepsilon^2},
j\leq k)$.
\smallskip
The following {\it snake property} is a consequence of the
definition of $W^\eps_s$: If $s<s'$ and $s$ and $s'$ belong to the same (open)
excursion interval of $\beta^\varepsilon$ away from $0$, then
$W^\eps_s(t)=W^\eps_{s'}(t)$
for every
$t\in [0,
\inf_{u\in[s, s']}\beta^\eps_u]$.

\medskip
\noindent{\bf 2.5 \ Convergence to super-Brownian motion}
\medskip
As in Section 1, we let $X^\varepsilon_t$ be $\varepsilon$ times the
sum of the point masses at the positions of the particles alive at time $t$
in the $\varepsilon$-system. This is equivalent to
writing
$$X^\varepsilon_t
=\int_0^{\tau_\varepsilon} dL^{\varepsilon,t}_s\,\delta_{W^\varepsilon_s(t)}.$$
To justify this formula, recall the correspondence between excursions and trees
described in Subsection 2.2 and note that each upcrossing time $s$
of $\beta^\varepsilon$ above level $t$ corresponds to one particle
alive at time $t$, whose position is $W^\varepsilon_s(t)$.
Similarly, the historical process $Y^\varepsilon_t$ is
$$Y^\varepsilon_t
=\int_0^{\tau_\varepsilon} dL^{\varepsilon,t}_s\,\delta_{W^\varepsilon_s}.$$

Recall our assumptions (1.1) and (1.2) and the convergence result in (1.3).
We next prove a result about
the uniform modulus of continuity for the paths $W^\varepsilon_s$.
For convenience, we make the convention that
$W^\varepsilon_s(t)=W^\varepsilon_s(\beta^\varepsilon_s)$ when
$t>\beta^\varepsilon_s$.

\medskip
\noindent{\bf Lemma 2.2}. {\it Let $\eta\in (0,{1\over 2})$. Then,
$$\lim_{\delta\downarrow 0}
\Big(\inf_{\varepsilon\in(0,1]}\,P\big[|W^\varepsilon_s(t+r)-W^\varepsilon_s(t)|
\leq
r^{{1\over 2}-\eta},
\ \hbox{for every }t\geq
0,\,r\in[0,\delta],\,s\in[0,\tau^\varepsilon]\big]\Big)=1.$$}

\noindent{\bf Remark.} This is of course reminiscent of the uniform modulus of
continuity for historical paths of super-Brownian motion. This lemma is
therefore very close to the results of [DIP] and [DP], which  however use
different approximations.
\smallskip
\noindent{\bf Proof.} Obviously it is enough to treat the case
when $x^\varepsilon_1=\cdots=x^\varepsilon_{N_\varepsilon}=0$ for every
$\varepsilon$.
We then use an embedding technique that will also play an important role
later. Let $(W_s,s\geq 0)$ be the Brownian snake of [L2] driven by the process
$(\beta_s,s\geq 0)$ and
with starting point $\underline{0}$. Recall that this is a continuous
Markov process with
values  in $\calw_0:=\{w\in\calw:w(0)=0\}$, whose law is characterized by
the following
properties:

\smallskip
$\bullet$ For every $s\geq 0$, the path $W_s$ has lifetime $\beta_s$.
\smallskip
$\bullet$ Conditionally on $(\beta_s,s\geq 0)$, the process $(W_s,s\geq 0)$
is time-inhomogeneous Markov, and its transition kernels are
characterized as follows. If $s<s'$, we have $W_{s'}(t)=W_s(t)$
for every $t\leq m(s,s'):=\inf_{[s,s']}\beta_r$, and
$(W_{s'}(m(s,s')+r)-W_{s'}(m(s,s')),0\leq r\leq \beta_{s'}-m(s,s'))$ is a
Brownian path
independent of $W_{s}$.
\smallskip
Now, for every $\varepsilon\in(0,1]$, we may assume that the spatial motions
of the particles are chosen in such a way that, for every $\eps>0$ and every
$k\in\{0,1,\ldots,
\tau^\eps/\eps^2\}$,
$$
\eqalign{W^\eps_{k\eps^2}&=W_{T^\eps_k}\hbox{ if } k \hbox{ is odd},\cr
W^\eps_{k\eps^2}&=W_{T^\eps_k}\mid [0,\beta_{T^\eps_k}-2\eps]\hbox{ if } k
\hbox{ is even},\cr}\eqno(2.6)
$$
where the notation $W_{T^\eps_k}\mid [0,\beta_{T^\eps_k}-2\eps]$ means that
the path $W_{T^\eps_k}$ is restricted to the interval
$[0,\beta_{T^\eps_k}-2\eps]
=[0,\beta^\eps_{k\eps^2}]$. In fact, it is immediate to verify that the
process $(W^\varepsilon_{k\varepsilon^2},0\leq k\leq
\tau^\varepsilon/\varepsilon^2)$
defined by (2.6)
has (conditionally on $\beta^\varepsilon$) the distribution described at the
end of Subsection 2.4.

Note that the family $(\tau^\varepsilon,\varepsilon\in(0,1])$
is bounded a.s.
Then the proof of Lemma 2.2 reduces to checking that, for every $K>0$,
$$\lim_{\delta\downarrow 0}
P\big[|W_s(t+r)-W_s(t)|\leq
r^{{1\over 2}-\eta},
\ \hbox{for every }t\geq 0,r\in[0,\delta],s\in[0,K]\big]=1.\eqno{(2.7)}$$
This can be easily done using Borel-Cantelli type arguments. Alternatively,
we may also use the relations between super-Brownian motion and the
Brownian snake [L2],
and the uniform modulus of continuity of [DP]. \hfill$\square$
\medskip
The graph ${\cal G}^\varepsilon$ of the $\varepsilon$-system of
branching particles is defined by
$${\cal G}^\varepsilon=
{\rm cl}\Big(\bigcup_{t\geq 0}\{t\}\times\supp X^\varepsilon_t\Big)
=\{W^\varepsilon_s(t): s\in[0,\tau^\varepsilon],0\leq t\leq
\beta^\varepsilon_s\}.
$$
We are interested in weak convergence of
${\cal G}^\varepsilon$ towards the graph ${\cal G}$ of $X$, which
we define as
$${\cal G}={\rm cl}\Big(\bigcup_{t\geq 0}\{t\}\times\supp X_t\Big).$$
We view both ${\cal G}^\varepsilon$ and ${\cal G}$ as random elements of
the space of all compact subsets of $\R_+\times \R$, which is
equipped with the Hausdorff metric.

\medskip
\noindent{\bf Lemma 2.3}. {\it We have the joint convergence
$$\big((X^\varepsilon_t,t\geq 0),{\cal G}^\varepsilon\big)\build\hbox to
10mm{\rightarrowfill} _{\varepsilon\to 0}^{\rm (d)}
\big((X_t,t\geq 0),{\cal G}).$$}
\smallskip
\noindent{\bf Proof.} We first consider the case
when $x^\varepsilon_1=\cdots=x^\varepsilon_{N_\varepsilon}=0$ for every
$\varepsilon$.
Then we can suppose that the processes
$(W^\varepsilon_s,s\in[0,\tau^\varepsilon])$ are
constructed via the embedding technique described in the
preceding proof. {}From (2.1) and (2.6), we get
$$(W^\varepsilon_{s\wedge \tau^\varepsilon},s\geq 0)
\build\hbox to
10mm{\rightarrowfill} _{\varepsilon\to 0}^{\rm (a.s.)}
(W_{s\wedge \tau},s\geq 0)\eqno{(2.8)}$$
in the sense of uniform convergence. Using
Lemma 2.1, we get
$$X^\varepsilon_t=\int_0^{\tau^\varepsilon}
dL^{\varepsilon,t}_s\,\delta_{W^\varepsilon_s(t)}
\build\hbox to
10mm{\rightarrowfill} _{\varepsilon\to 0}^{\rm (a.s.)}
\int_0^\tau dL^t_s\,\delta_{W_s(t)}=X_t$$
uniformly in $t$. (The formula for $X_t$ is the Brownian snake
representation of super-Brownian motion, see [L2].) Furthermore, the
convergence
(2.8) also implies that
$${\cal G}^\varepsilon=\{W^\varepsilon_s(t):s\leq \tau^\varepsilon,
t\leq \beta^\varepsilon_s\}\build\hbox to
10mm{\rightarrowfill} _{\varepsilon\to 0}^{\rm (a.s.)}
\{W_s(t):s\leq \tau,
t\leq \beta_s\},$$
and the limit is easily identified with the graph ${\cal G}$ of $X$.
Therefore we get the statement of the lemma in the special
case $x^\varepsilon_1=\cdots=x^\varepsilon_{N_\varepsilon}=0$.

Before proceeding to the general case, let us make one more
observation. Fix $\delta>0$ and write $(V^\varepsilon_s,s\geq 0)$
for a process distributed as an excursion of $W^\varepsilon$
away from $\underline 0$ conditioned to have height greater than $\delta$.
(Alternatively, $(V^\varepsilon_s,s\geq 0)$ codes the historical paths
of the $\varepsilon$-system starting with one particle at the
origin and conditioned to be non-extinct at time $\delta$.)
It follows from the convergence (2.8) that we have also
$$(V^\varepsilon_{s},s\geq 0)
\build\hbox to
10mm{\rightarrowfill} _{\varepsilon\to 0}^{\rm (d)}
(V_{s},s\geq 0),$$
where the limiting process is an excursion
of $W$ conditioned to have height greater than $\delta$.
As in the first part of the proof, it follows that the graphs of
$V^\varepsilon$
(defined analogously to ${\cal G}^\varepsilon$) also converge in
distribution towards the graph of $V$. Furthermore, this convergence holds
jointly
with that of the measure-valued processes ${\cal X}^\varepsilon_t$
associated with
$V^\varepsilon$ in the same way as $X^\varepsilon_t$ was associated with
$W^\varepsilon$.

Let us consider now the general case. Because of Lemma 2.2 and
assumption (1.2), it is enough to
prove that for any fixed $\delta>0$, ${\cal
G}_\varepsilon\cap([\delta,\infty)\times\R)$
converges in distribution to ${\cal G}\cap([\delta,\infty)\times\R)$
(and that this convergence holds jointly with that of $X^\varepsilon$).
Let $A_\varepsilon$ stand for the set of indices
$j\in\{1,\ldots,N_\varepsilon\}$
such that the $j$-th excursion of $\beta^\varepsilon$ has a height greater
than
$\delta$. Note that the events $\{j\in A_\varepsilon\}$ are independent
with the
same probability $2\varepsilon(2\varepsilon+\delta)^{-1}$. It follows that the
random measure
$$\sum_{j\in A_\varepsilon} \delta_{x^\varepsilon_j}$$
converges weakly to a Poisson measure with intensity ${2\over \delta}\mu$.
Note that,
conditionally on $A_\varepsilon$, ${\cal
G}_\varepsilon\cap[\delta,\infty)\times\R$
has the same distribution as
$$\bigcup_{j\in A_\varepsilon} \big(((0,x^\varepsilon_j)+{\cal
G}_{j,\varepsilon})
\cap[\delta,\infty)\times\R\big)$$
where ${\cal G}_{j,\varepsilon}$ are independent copies of the graph of
$V^\varepsilon$. If follows that the random sets
${\cal G}_\varepsilon\cap[\delta,\infty)\times\R$ converge in distribution to
$$\bigcup_{j\in J} \big(((0,x_j)+{\cal G}_{(j)})
\cap[\delta,\infty)\times\R\big),$$
where $\sum_{j\in J}\delta_{x_j}$ is a Poisson point measure on $\R$
with intensity ${2\over \delta}\mu$, and, conditionally on this random measure,
the random sets ${\cal G}_{(j)}$ are independent and distributed according to
the law of the graph of $V$. The canonical representation of
superprocesses allows us to identify this limiting distribution with that
of ${\cal G}\cap ([\delta,\infty)\times\R)$. Furthermore,
using the joint convergence of $(V^\varepsilon,{\cal X}^\varepsilon)$,
it is easy to verify
that the convergence holds jointly with that of $X_\varepsilon$.
\hfill$\square$

\medskip
\noi{\bf 3. Branching particle systems with reflection}
\medskip
\noi{\bf 3.1 \ Reflection for deterministic paths}
\medskip
The purpose of this section is to explain, first in a deterministic
setting, the construction of reflected systems.
We consider a deterministic branching particle system in $\R$ analogous
to the ones considered above. At time $0$,
we have $N$ particles located at $x_1,\ldots,x_N$.
Each particle moves in $\R$ and gives birth at its death
to $0$ or $2$ new particles. As in Subsection 2.2, denote by
${\cal T}$ the genealogical forest
of the population,
which is a subset of $\{1,\ldots,N\}\times {\bf U}$. Each element $v=(k,u)$
in ${\cal T}$ corresponds to a particle with birth time $\xi_v$
and death time $\zeta_v$
(as in Section 2, we could alternatively consider the life durations
$\ell_v:=\zeta_v-\xi_v$ but in this subsection and the
next one it is more convenient to deal with the birth and
death times). The spatial motion of $v$ is a
continuous function $f_v:[\xi_v,\zeta_v]\longrightarrow \R$
and $f_{v'}(\xi_{v'})=f_v(\zeta_v)$ if $v'$ is a child of $v$
(then $\xi_{v'}=\zeta_v$). The historical path of $v$
is the continuous function $w_v:[0,\zeta_v]\longrightarrow \R$
such that, for every $t\in[0,\zeta_v)$, $w_v(t)$ is the position
at time $t$ of the ancestor of $v$ alive at that time.

We
assume that the death times $\zeta_v$, $v\in{\cal T}$ are all distinct, that
the system becomes extinct after a finite number
of generations and that when a particle dies there
is no other particle at the same location: For every $v\in {\cal T}$,
$f_v(\zeta_v)\not =f_{v'}(\zeta_v)$ for every $v'\in{\cal T}$ such
that $\xi_{v'}\leq \zeta_v<\zeta_{v'}$.

We turn to the construction of the reflected system. This system is such that
the number and positions of the particles alive at every time $t$ are the
same as
in the original system (thus each death time for the reflected system
is also a death time for the reflected system). However the genealogical
forest $\wt{\cal T}$ will be different, as will be the spatial
motions $\wt f_u,\,u\in\wt{\cal T}$ or the birth and death times
$\wt\xi_u,\,\wt\zeta_u,\,u\in\wt{\cal T}$.

Set $R_0=0$ and denote by $R_1<R_2<\cdots<R_M$ the successive death times
in the original
system.
For every $k\in\{1,\ldots,M\}$, let ${\cal T}_{(k)}$ be the set of (labels
of) particles
that are alive on the interval $[R_{k-1},R_k)$.
We use induction on $k$ to define sets $\wt{\cal T}_{(k)}$, which will
represent the particles alive on the interval $[R_{k-1},R_k)$ in the
reflected system, and the corresponding spatial motions.

To begin with, we have $\wt{\cal T}_{(1)}=\{1,\ldots,N\}$, and
we define $\wt f_j(t)$ for every $t\in [0,R_1]$ and every $j\in\wt{\cal
T}_{(1)}$
by requiring $(\wt f_1(t),\ldots,\wt f_N(t))$ to be the
increasing rearrangement of $(f_1(t),\ldots,f_N(t))$. Note that the mappings
$\wt f_1,\ldots,\wt f_N$ are continuous.

Suppose that for some $k\in\{1,\ldots, M-1\}$, we have defined
$\wt{\cal T}_{(k)}$ and the corresponding paths $(\wt
f_u(t),t\in[R_{k-1},R_k])$,
for $u\in\wt{\cal T}_{(k)}$, in such a way that $\card\,\wt{\cal T}_{(k)}=
\card\,{\cal T}_{(k)}$, and, for
every  $t\in[R_{k-1},R_k]$:
\smallskip
$\bullet$ The
mapping $\wt{\cal T}_{(k)}\ni u\to \wt f_u(t)$
is increasing with respect to the lexicographical order on $\wt{\cal T}_{(k)}$.
\smallskip
$\bullet$ The
values of $\wt f_u(t)$ for
$u\in
\wt{\cal T}_{(k)}$ (counted with their multiplicities)
are the same as those of $f_u(t)$ for $u\in {\cal T}_{(k)}$.
\smallskip
By definition, one of the particles in ${\cal T}_{(k)}$,
say $u_{(k)}$, dies at time $R_k$. Then there is exactly one
$\wt u_{(k)}\in\wt{\cal T}_{(k)}$ such that $\wt f_{\tilde
u_{(k)}}(R_k)=f_{u_{(k)}}(R_k)$.
We set
$$\wt{\cal T}_{(k+1)}=\big(\wt{\cal T}_{(k)}\backslash \{\wt u_{(k)}\}\big)
\cup \{\wt u_{(k)}1,\wt u_{(k)}2\}$$
if $u_{(k)}$ has two children in the original system, and
$$\wt{\cal T}_{(k+1)}=\wt{\cal T}_{(k)}\backslash \{\wt u_{(k)}\}$$
if not. Furthermore, let $u^{k+1}_1,\ldots,u^{k+1}_{N_{k+1}}$
be the elements of $\wt{\cal T}_{(k+1)}$ listed in lexicographical order.
We define $\wt f_u(t)$ for every $t\in [R_k,R_{k+1}]$ and every
$u\in\wt{\cal T}_{(k+1)}$
by requiring that $(\wt f_{u^{k+1}_1}(t),\ldots,\wt
f_{u^{k+1}_{N_{k+1}}}(t))$ is the
increasing rearrangement of $(f_u(t),u\in{\cal T}_{(k+1)})$. Notice that
when $u\in \wt{\cal T}_{(k)}\cap \wt{\cal T}_{(k+1)}$
the definition of $\wt f_u(R_k)$ is consistent with the previous step.

Finally, the genealogical forest of the reflected system is
$$\wt{\cal T}=\bigcup_{k=1}^M \wt{\cal T}_{(k)}.$$
The birth and death times $\wt\xi_u,\wt\zeta_u$
as well as the (continuous) spatial motions $\wt f_u$
in the reflected system are defined by the requirement of consistency with
the construction of $\wt{\cal T}_{(k)}$'s. Note the two fundamental properties:
\smallskip
$\bullet$ At each time $t\geq 0$, the positions of the particles
(counted with their multiplicities) are the same in the original
and the reflected system.
\smallskip
$\bullet$ If $u,v\in\wt{\cal T}$
 with $u\prec v$ ($\prec$ denotes the lexicographical order) then
$f_u(t)\leq f_v(t)$ for every $t\in[\wt\xi_u,\wt\zeta_u]\cap
[\wt\xi_v,\wt\zeta_v]$.
\smallskip
Historical paths $\wt w_u$, $u\in\wt{\cal T}$
for the reflected system
are defined in a way analogous to the original one.
If $u,v\in\wt{\cal T}$ and $u\prec v$ then $\wt w_u(t)\leq \wt w_{v}(t)$
for every
$t\in[0,\wt\zeta_u\wedge\wt\zeta_v]$.
\medskip
\noi{\bf 3.2 \ A technical lemma}
\medskip
Let $M\in\{1,\cdots,N\}$, and consider a branching system
consisting only of the particles
labeled $1,\ldots,M$ at time $0$ and their descendants. The new genealogy
is described by the forest
$${\cal T}':={\cal T}\cap(\{1,\ldots,M\}\times U).$$
{}From this new branching particle system, we can construct a reflected system
by the procedure described in Subsection 3.1. We denote by $\wt{\cal T}'$ the
genealogical forest for this new reflected system, and by $\wt w'_v$, $v\in
\wt{\cal T}'$
the associated historical paths. In general, the historical paths
$\wt w'_v$ will be very different from those obtained by reflecting the
original
system. Under special assumptions however, we can say that some of the
paths $\wt w'_v$ will also be (reflected)
historical paths in the original system.
\medskip
\noi{\bf Lemma 3.1}. {\it Let $t>0$ and let $I$ be a bounded interval in
$\R$. Suppose
that $w_v(r)\notin I$ for every $v\in{\cal T}\backslash{\cal T}'$
and $r\in[0,t]$. If $v\in\wt{\cal T}$ is such that $\wt\zeta_v\geq t$
and $\wt w_v(r)\in I$ for every $r\in[0,t]$,
then there exists
$v'\in\wt{\cal T}'$ such that $\wt\zeta'_{v'}\geq t$ and $\wt
w'_{v'}(r)=\wt w_v(r)$
for every $r\in[0,t]$. The converse also holds: If
$v'\in\wt{\cal T}'$ is such that $\wt\zeta'_{v'}\geq t$ and $\wt
w'_{v'}(r)\in I$ for every
$r\in[0,t]$, then there exists $v\in\wt{\cal T}$ such that $\wt\zeta_v\geq t$
and $\wt w_v(r)=\wt w'_{v'}(r)$
for every $r\in[0,t]$.}
\medskip
In other words, the first assertion means that the path $\wt w_v$, or
rather its
restriction to $[0,t]$, will still be a historical path for the new reflected
system. We leave an easy proof of the lemma to the reader.
\medskip
\noi{\bf 3.3 \ Reflected branching particle systems}
\medskip
For every $\varepsilon\in(0,1]$,
we can apply the construction of Subsection 3.1 to the
$\varepsilon$-system of branching Brownian motions. Note that the
assumptions that we
imposed on the deterministic system hold with probability one for this
random system. We
write
${\cal T}_\varepsilon$ for the genealogical forest of the
$\varepsilon$-system, and $(\ell^\varepsilon_u,u\in{\cal T}_\varepsilon)$
for the lifetimes of particles. The notation
$\wt{\cal T}_\varepsilon$ and $(\wt\ell^\varepsilon_u,u\in\wt{\cal
T}_\varepsilon)$
has a similar meaning for the corresponding reflected system,
which we call the $\varepsilon$-reflected system. Observe that $({\cal
T}_\varepsilon,(\ell^\varepsilon_u,u\in{\cal T}_\varepsilon))$ and
$(\wt{\cal T}_\varepsilon,(\wt\ell^\varepsilon_u,u\in\wt{\cal
T}_\varepsilon))$ have the
same distribution. This is so because the spatial
motions and branching structure for the $\varepsilon$-system
of branching Brownian motions are independent (a tedious rigorous
justification could be given, but we feel that the result is
sufficiently obvious to allow us to omit it). Furthermore,
$\card\,\wt{\cal T}_\varepsilon=\card\,{\cal T}_\varepsilon$.
\smallskip
We noticed at the end of Subsection 2.2 that the
process $(\beta^\varepsilon_s,s\in[0,\tau^\varepsilon])$
can be reconstructed as a measurable function of the marked trees $({\cal
T}_\varepsilon,(\ell^\varepsilon_u,u\in{\cal T}_\varepsilon))$. Hence, we can
also code the
branching structure of the $\varepsilon$-reflected system by
a random process
$(\wt\beta^\varepsilon_s,s\in[0,\wt\tau^\varepsilon])$ which has the
same distribution as $(\beta^\varepsilon_s,s\in[0,\tau^\varepsilon])$.
The fact that
$\card\,\wt{\cal T}_\varepsilon=\card\,{\cal T}_\varepsilon$
implies that the time $\tau_\varepsilon$ is also the end of
the $N_\varepsilon$-th excursion of $\wt\beta^\varepsilon$
away from $0$, and thus $\wt\tau^\varepsilon=\tau^\varepsilon$. The
discrete local
times  of $\wt\beta^\varepsilon$ (cf. Subsection 2.3) are
denoted by $(\wt L^{\varepsilon,x}_s,x\in\R_+,s\in[0,\tau^\varepsilon])$.
\smallskip
Finally, we can code the historical paths of the $\varepsilon$-reflected system
by a discrete snake
$(\wt W^\eps_s, s\in[0, \tau^\eps])$ in a way analogous to what we did in
Subsection 2.4.
Recall that we assume $x^\varepsilon_1\leq \cdots\leq
x^\varepsilon_{N_\varepsilon}$.
As in Section 2, if $s\in\eps^2\N\cap[0,\tau^\varepsilon)$ and
$\wt\beta^\varepsilon_s=0$,
we set $\wt W^\varepsilon_s=\underline x^\varepsilon_k$ if $s$ is the
beginning of the $k$-th excursion of $\wt\beta^\varepsilon$ away from $0$
(and $\wt W^\varepsilon_{\tau_\varepsilon}=\underline
x^\varepsilon_{N_\varepsilon}$).
Otherwise, if $s\in\eps^2\N\cap[0,\tau^\varepsilon)$ and
$\wt\beta^\eps_s>0$, then
$(s, \wt\beta^\eps_s)$ can be associated
with a unique edge $u$ of the forest $\wt{\cal T}_\varepsilon$, and
we let $\wt W^\varepsilon_s$ be equal to $\wt w^\varepsilon_u$,
the historical path of $u$. If
$s\not\in\eps^2\N$, we use the same interpolation as in Section 2.
A fundamentally important property of the
process $(\wt W^\eps_s, s\in[0, \tau^\eps])$, from the point of view
of our project, is that for $s<s'$,
$$\wt W^\eps_s(t)\le\wt W^\eps_{s'}(t), \qquad \forall\;t\in[0,
\wt\beta^\eps_s\wedge\wt\beta^\eps_{s'}].\eqno{(3.1)}$$
This follows from our construction and the end of Subsection 3.1.
As in the case of $W^\eps_s$, we see that
if $s<s'$, then
$$\wt W^\eps_s(t)=\wt W^\eps_{s'}(t),
\qquad\forall \
t\in [0, \inf_{u\in[s, s']}\wt \beta^\eps_u].$$

Because at every time the locations of particles
are the same in the reflected system and in the original one,
the random measure
$$\wt X^\varepsilon_t=\int_0^{\tau_\varepsilon} d\wt
L^{\varepsilon,t}_s\,\delta_{
\wt W^\varepsilon_s(t)}$$
coincides with $X^\varepsilon_t$. Things are however very different for the
historical measure
$$\wt Y^\varepsilon_t=\int_0^{\tau_\varepsilon} d\wt
L^{\varepsilon,t}_s\,\delta_{
\wt W^\varepsilon_s}.$$

\medskip
\noi{\bf 4. Tightness of the reflected system}
\medskip
\noi{\bf 4.1 \ Uniform continuity of the reflected paths}
\medskip
Our first goal is to derive an important uniform
continuity property for the individual
paths of the $\varepsilon$-reflected
system (Theorem 4.1 below).
{}From the intuitive point of view, reflected paths
should have smaller oscillations than ``free'' paths
and so this property seems to be a straightforward
consequence of Lemma 2.2.
However the intuition about the relationship between
moduli of continuity of free and reflected paths
is only  correct as long as we do not have any deaths.
To be specific consider $p$
paths $w_{(1)},\ldots,w_{(p)}$
all defined on the time interval $[0,1]$, and let
$\wt w_{(1)},\ldots,\wt w_{(p)}$ be the corresponding system
of reflected paths. Then, if we assume that
$|w_{(i)}(t)-w_{(i)}(t')|\leq \varphi(|t-t'|)$
for every $i=1,\ldots,p$ and $t,t'\in[0,1]$ and for some
nondecreasing function $\varphi$, an easy argument shows that the same
bound holds
when the paths $w_{(i)}$ are replaced by $\wt w_{(i)}$.

It turns out that a similar assertion about moduli of continuity
is false if paths may have different lifetimes.
Fig.~2 shows a system of two paths. In the
original system, the oscillations of paths
over the intervals where they are defined
are equal to $z_1 - y_1$ and $z_2 - y_2$.
One of the paths in the reflected system
goes from $y_1$ to $z_2$ and so has an oscillation
larger than the oscillations of the original paths.
In this article we consider Brownian particles which die
at different times so we cannot use known estimates
for the modulus of continuity of the original
(non-reflecting) historical paths in a direct way.
We will use them later in a different but quite essential way.

\bigskip
\vbox{
\epsfxsize=2.5truein
  \centerline{\epsffile{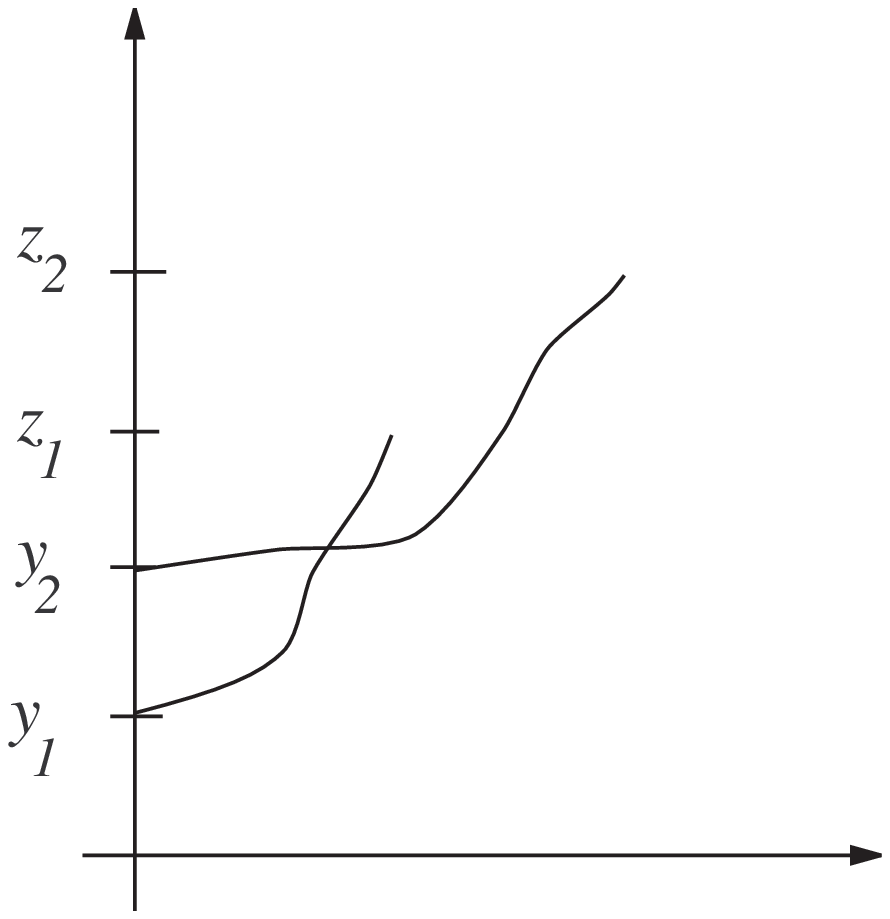}}

\centerline{Figure 2.}
}
\bigskip

Recall our notation $\wt w^\varepsilon_u$, $u\in\wt{\cal T}_\varepsilon$,
for the historical paths of the $\varepsilon$-reflected system. By convention,
$\wt w^\varepsilon_u(t)=\wt w^\varepsilon_u(\wt\zeta^\varepsilon_u)$
if $t\geq \wt\zeta^\varepsilon_u$.
\smallskip
\noi{\bf Theorem 4.1}. {\it For every $\eta>0$,
$$\lim_{\delta\to 0}\Big(
\limsup_{\varepsilon\to 0}\,P\Big[\build{\sup_{t,t'\geq 0}}_{|t-t'|\leq
\delta}^{}
\,\sup_{u\in\wt{\cal T}_\varepsilon}\,|\wt w^\varepsilon_u(t)-
\wt w^\varepsilon_u(t')|>\eta\Big]\Big)=0.$$}
\smallskip
\noi{\bf Proof.} Let $(\delta_{(p)},\varepsilon_{(p)})$ be a sequence in
$(0,1]^2$
converging to $0$. We will prove that there exists a subsequence
$(\delta'_{(p)},\varepsilon'_{(p)})$ such that:
$$\lim_{p\to\infty}\Big(\,
\build{\sup_{t,t'\geq 0}}_{|t-t'|\leq \delta'_{(p)}}^{}
\,\sup_{u\in\wt{\cal T}_{\varepsilon'_{(p)}}}\,|\wt
w^{\varepsilon'_{(p)}}_u(t)-
\wt w^{\varepsilon'_{(p)}}_u(t')|\Big)=0\eqno{(4.1)}$$
in probability. Clearly, the statement of Theorem 4.1 is a consequence
of this fact.

We first explain how we choose the sequence
$(\delta'_{(p)},\varepsilon'_{(p)})$.
By Lemma 2.3 and the Skorohod representation theorem ([EK] Theorem 3.1.8),
we may,
for every $p\geq 1$, replace
the pair $(X^{\varepsilon_{(p)}},{\cal G}_{\varepsilon_{(p)}})$
by a new pair with the same distribution (for which we keep the same
notation), in such
a way that
$$(X^{\varepsilon_{(p)}},{\cal G}_{\varepsilon_{(p)}})
\build\hbox to
10mm{\rightarrowfill} _{p\to\infty}^{\rm (a.s.)}
(X,{\cal G}),$$
where $X$ is a super-Brownian motion started at $\mu$ and ${\cal G}$
denotes its graph.
Note that the genealogical forest $\wt{\cal T}_{\varepsilon_{(p)}}$, the
process $\wt \beta^{\varepsilon_{(p)}}$
and the historical paths $\wt w^{\varepsilon_{(p)}}_u$, $u\in\wt{\cal
T}_{\varepsilon_{(p)}}$, are reconstructed as
measurable functions of the new process $X^{(\varepsilon_p)}$, and that it
suffices to prove (4.1) for the new historical paths. As a consequence of the
remark following Lemma 2.1, we have
$$\lim_{p\to\infty}\Big(\;\sup_{s\geq 0}
\build{\sup_{t,t'\geq 0}}_{|t-t'|\leq \delta_{(p)}}^{}|\wt
L^{\varepsilon_{(p)},t}_{s\wedge
\tau^{\varepsilon_{(p)}}}-\wt L^{\varepsilon_{(p)},t'}_{s\wedge
\tau^{\varepsilon_{(p)}}}|
\Big)=0\eqno{(4.2)}$$
in probability. We choose the subsequence $(\delta'_{(p)},\varepsilon'_{(p)})$
so that the convergence (4.2) holds almost surely along this subsequence.

We will argue by contradiction to prove (4.1). If (4.1) does not hold,
then on a set $A$ of positive probability, we can find a number $\eta>0$
and a (random)
subsequence $p_k\uparrow \infty$ such that, if
$\varepsilon_k:=\varepsilon'_{(p_k)}$
and $\delta_k:=\delta'_{(p_k)}$,
$$\build{\sup_{t,t'\geq 0}}_{|t-t'|\leq \delta_{k}}^{}
\,\sup_{u\in\wt{\cal T}_{\varepsilon_{k}}}\,|\wt w^{\varepsilon_{k}}_u(t)-
\wt w^{\varepsilon_{k}}_u(t')|>\eta.\eqno{(4.3)}$$
{}From now on until the end of the proof, we will assume that
the event $A$ holds.
By (4.3), for every $k\geq 1$, there exist $u_k\in\wt{\cal T}_{\varepsilon_k}$,
$t_k,t'_k\geq 0$ with $|t_k-t'_k|\leq \delta_k$, such that
$$|\wt w^{\varepsilon_k}_{u_k}(t_k)-\wt w^{\varepsilon_k}_{u_k}(t'_k)|>\eta.$$
Clearly, we can assume that $t_k\leq
t'_k\leq \wt\zeta^{\varepsilon_k}_{u_k}$.

Recall that the graphs ${\cal G}_{\varepsilon_{(p)}}$ converge to
${\cal G}$
in the Hausdorff metric. In particular, the set of all pairs $(t_k,\wt
w^{\varepsilon_k}
_{u_k}(t_k))$ and $(t'_k,\wt w^{\varepsilon_k}
_{u_k}(t'_k))$ is relatively compact.
By passing to a subsequence, if necessary, we
may assume that $t_k,t'_k\longrightarrow t_\infty$, $\wt w^{\varepsilon_k}
_{u_k}(t_k)\longrightarrow x_1$ and $\wt w^{\varepsilon_k}
_{u_k}(t'_k)\longrightarrow x_2$ as $k\to \infty$. We have $|x_1-x_2|\geq
\eta$,
and we take $x_2>x_1$ for definiteness.

We also know that $\wt X^{\varepsilon_{(p)}}_t=X^{\varepsilon_{(p)}}_t$
converges to $X_t$ a.s. as $p\to\infty$, uniformly on compact subsets
of $\R_+$. Hence, both sequences $\wt X^{\varepsilon_k}_{t_k}$ and
$\wt X^{\varepsilon_k}_{t'_k}$ converge to $X_{t_\infty}$, and
$$\eqalign{
&\limsup_{k\to\infty} \wt X^{\varepsilon_k}_{t_k}((-\infty,\wt
w^{\varepsilon_k}
_{u_k}(t_k)])\leq X_{t_\infty}((-\infty,x_1])\,,\cr
&\liminf_{k\to\infty} \wt X^{\varepsilon_k}_{t'_k}((-\infty,\wt
w^{\varepsilon_k}
_{u_k}(t'_k)))\geq X_{t_\infty}((-\infty,x_2))\,.}\eqno{(4.4)}
$$
We claim that
$$\liminf_{k\to\infty}
\Big(\wt X^{\varepsilon_k}_{t_k}\big((-\infty,\wt w^{\varepsilon_k}
_{u_k}(t_k)]\big)-\wt X^{\varepsilon_k}_{t'_k}\big((-\infty,\wt
w^{\varepsilon_k}
_{u_k}(t'_k))\big)\Big)\geq 0.\eqno{(4.5)}$$
To see this, we use the discrete snake representation
of Subsection 3.3. Write $s_k\in [0,\tau_{\varepsilon_k})\cap \varepsilon^2\N$
for the time associated with the edge $u_k$
of $\wt {\cal T}_{\varepsilon_k}$ in this representation. By construction, $\wt
w^{\varepsilon_k} _{u_k}=\wt W^{\varepsilon_k}_{s_k}$, and (3.1) implies
$$\wt X^{\varepsilon_k}_{t_k}((-\infty,\wt w^{\varepsilon_k}
_{u_k}(t_k)])=\int_0^{\tau_{\varepsilon_k}}
d\wt L^{\varepsilon_k,t_k}_s\,
{\bf 1}_{\{\wt W^{\varepsilon_k}_{s}(t_k)
\leq \wt W^{\varepsilon_k}_{s_k}(t_k)\}}\geq \wt L^{\varepsilon_k,t_k}_{s_k}.$$
Similarly, we get
$$\wt X^{\varepsilon_k}_{t'_k}((-\infty,\wt w^{\varepsilon_k}
_{u_k}(t'_k)))\leq \wt L^{\varepsilon_k,t'_k}_{s_k}.$$
Hence,
$$\wt X^{\varepsilon_k}_{t_k}((-\infty,\wt w^{\varepsilon_k}
_{u_k}(t_k)])-\wt X^{\varepsilon_k}_{t'_k}((-\infty,\wt w^{\varepsilon_k}
_{u_k}(t'_k)))\geq \wt L^{\varepsilon_k,t_k}_{s_k}-\wt
L^{\varepsilon_k,t'_k}_{s_k}.
\eqno{(4.6)}$$
On the other hand,
$$|\wt L^{\varepsilon_k,t_k}_{s_k}-\wt L^{\varepsilon_k,t'_k}_{s_k}|
\leq \sup_{s\in[0,\tau_{\varepsilon_k}]}\,
\build{\sup_{t,t'\geq 0}}_{|t-t'|\leq \delta_k}^{}|\wt L^{\varepsilon_k,t}_s-
\wt L^{\varepsilon_k,t'}_s|.$$
By the convergence in (4.2), which holds a.s. along
the subsequence $(\delta'_{(p)},\varepsilon'_{(p)})$, the right hand side tends
to $0$
as $k\to\infty$. This and (4.6) give the claim (4.5).

{}From (4.5) and (4.4), we get $X_{t_\infty}((-\infty,x_1])\geq
X_{t_\infty}((-\infty,x_2))$
and thus (recall that $x_1<x_2$), $X_{t_\infty}((x_1,x_2))=0$.
This a priori does not imply that $\{t_\infty\}\times (x_1,x_2)\cap {\cal
G}=\emptyset$
as there could be a ``local extinction'' of $X$ at time $t_\infty$ in
$(x_1,x_2)$.
However, by Theorem 1.4 of Perkins [P1], there can be at most
one local extinction at a given time, so we can choose
$x'_1$ and $x'_2$ with $x_1<x'_1<x'_2<x_2$ such that
$\{t_\infty\}\times [x'_1,x'_2]\cap {\cal G}=\emptyset$. Since ${\cal G}$
is closed, we have also $[t_\infty-\delta,t_\infty+\delta]\times
[x'_1,x'_2]\cap {\cal
G}=\emptyset$ for $\delta>0$ sufficiently small. However, by construction,
for $k$ sufficiently large the paths $\wt w^{\varepsilon_k}_{u_k}$
and thus also the graph ${\cal G}_{\varepsilon_k}$
must intersect $[t_\infty-\delta,t_\infty+\delta]\times [x'_1,x'_2]$. This
gives
a contradiction since we know that ${\cal G}_{\varepsilon_k}$
converge to ${\cal G}$. This contradiction completes
the proof of Theorem 4.1. \hfill$\square$

\medskip
\noi{\bf 4.2 \ Tightness of reflected discrete snakes}
\medskip
{}From now on, we restrict our attention to
values of $\varepsilon$ belonging to a fixed sequence
$\cE$ decreasing to $0$. For convenience, we extend the
definition of the discrete snakes $\wt W^\varepsilon$
by taking $\wt W^\varepsilon_s=\wt
W^\varepsilon_{\tau^\varepsilon}=\underline{x}
^\varepsilon_{N_\varepsilon}$ (and thus $\wt\beta^\varepsilon_s=0$) for
$s>\tau^\varepsilon$.
\medskip
\noi {\bf Proposition 4.2}. {\it The laws of the processes $\wt W^\varepsilon$,
$\varepsilon\in \cE$, are tight in the space of all
probability measures on $\D([0,\infty),{\cal W})$. Furthermore,
if $(\wt W_s,s\geq 0)$ is a weak limit point of this sequence of processes,
we have the following properties.

{\rm (i)} If $\wt \beta_s:=\zeta_{\wt W_s}$, the process $(\wt
\beta_s,s\geq 0)$
has the same distribution as $(\beta_{s\wedge \tau},s\geq 0)$.

{\rm (ii)} Almost surely for every $s\leq s'$ we have $\wt W_s(t)\leq \wt
W_{s'}(t)$
for every $t\in[0,\wt\beta_s\wedge\wt\beta_{s'}]$.

{\rm (iii)} The set of discontinuities of the mapping $s\to \wt W_s$
is contained in the zero set of $\wt \beta$. Furthermore, if
$s<s'$ belong to the same connected component of the complement of
the zero set, we have
$$\wt W_s(t)=\wt W_{s'}(t)\qquad\hbox{for every }
t\in[0,\inf_{r\in[s,s']}\wt\beta_r]\,.$$}
\noi{\bf Proof}. The hard part of the proof is to show tightness.
To this end we rely
on the classical criteria (see e.g. Corollary 3.7.4 of [EK]). We first observe
that the compact containment condition is a straightforward
consequence of Theorem 4.1. In fact, if $\eta>0$ is fixed,
then for every integer $p\geq 1$, Theorem 4.1
and the construction of the discrete snake $\wt W^\varepsilon$
allow us to find $\delta_p>0$ such that, for $\varepsilon\in \cE$
small enough,
$$P\Big[\sup_{s\geq 0}\,\build{\sup_{t,t'\geq 0}}_{|t-t'\leq \delta_p}^{}
|\wt W^\varepsilon_s(t)-\wt W^\varepsilon_s(t')|>2^{-p}\Big]\leq
\eta\,2^{-p-1}.
\eqno{(4.7)}$$
(Here and later, we make the convention that $\wt W^\varepsilon_s(t)
=\wt W^\varepsilon_s(\wt\beta^\varepsilon_s)$ for $t>\wt\beta^\varepsilon_s$.)
It is easy to see that an even stronger assertion holds, namely,
(4.7) is true for all $\varepsilon\in \cE$; this can be achieved
by taking $\delta_p$
even smaller if necessary---note that for any fixed value of $\varepsilon$
we need only consider a finite number of historical paths. Then let
$H$ be a compact subset of $\R_+$
containing $\supp \mu_\varepsilon$ for $\varepsilon\in\cE$, and let $A>0$
be a constant. The set
$$\eqalign{K:=&\{w\in{\cal W}: w(0)\in H\,,\zeta_w\leq A\,,\cr
&\ \hbox{and }|w(t)-w(t')|\leq 2^{-p}\hbox{ for every }t,t'\in[0,\zeta_w]
\hbox{ with }|t-t'|\leq \delta_p\hbox{ and every }p\geq 1\}}$$
is compact, and it follows from (4.7) that
$$P[\wt W^\varepsilon_s\notin K \hbox{ for some }s\geq 0]<\eta$$
provided that $A$ is chosen large enough.

Recall the definition of the distance $d$ from Subsection 2.4. We set
$$\theta(\varepsilon,\delta)
=\inf_{(s_i)}\Big\{\sup_{i}\,\sup_{s,s'\in[s_{i-1},s_i)}
d(W^\varepsilon_s,W^\varepsilon_{s'})\Big\}\,,$$
where the infimum is over all finite sequences
$0=s_0<s_1<\cdots<s_{m-1}<\tau^\varepsilon\leq s_m$
such that $\inf\{|s_i-s_{i-1}|;1\leq i\leq m\}\geq \delta\}$.
As a direct application of Corollary 3.7.4 in [EK], the proof
of tightness will be complete if we can verify that, for every
$\eta>0$, we can choose $\delta>0$ sufficiently small
so that
$$\limsup_{\cE\ni\varepsilon\to 0}
P[\theta(\varepsilon,\delta)>\eta]<\eta.
\eqno{(4.8)}$$

We now fix $\eta>0$ and proceed to the proof of (4.8). As a consequence of
Theorem 4.1,
we can choose $\rho\in(0,\eta/5)$ so small that, for every $\varepsilon\in\cE$,
$$P\Big[\sup_{s\geq 0}\,\build{\sup_{t,t'\geq 0}}_{|t-t'|\leq \rho}^{}
|\wt W^\varepsilon_s(t)-\wt W^\varepsilon_s(t')|\leq {\eta\over 5}\Big]\geq
1-{\eta\over
5}.
\eqno{(4.9)}$$
Then, by the tightness of the laws of $\wt \beta^\varepsilon$
(cf (2.3)), we can choose $\kappa>0$
small enough so that, for every $\varepsilon\in\cE$,
$$P\Big[\build{\sup_{s,s'\geq 0}}_{|s-s'|\leq \kappa}^{}
|\wt \beta^\varepsilon_s-\wt \beta^\varepsilon_{s'}|\leq \rho\Big]\geq
1-{\eta\over 5}.
\eqno{(4.10)}$$
We denote by $E_\varepsilon$ the intersection of the events considered in
(4.9) and (4.10), so that the probability of the complement of
$E_\varepsilon$ is bounded above by $2\eta/5$.

Set $\gamma=\eta/5$. Since $\mu$ is a finite measure
with compact support, we can easily find an integer $M_\gamma$
and a finite sequence of reals $y_1<z_1\leq y_2<z_2\leq\cdots\leq y_{M_\gamma}
<z_{M_\gamma}$, such that:

$\bullet$ $z_i-y_i<\gamma$ for every $i=1,\ldots , M_\gamma$,

$\bullet$ $\displaystyle{\bigcup_{i=1}^{M_\gamma}[y_i,z_i)}$
contains a neighborhood of $\supp \mu$,

$\bullet$ $\mu(\{y_i\})=\mu(\{z_i\})=0$, and $\mu([y_i,z_i))>0$
for every $i=1,\ldots M_\gamma$.

\smallskip

\noi By the last condition,
$a_\gamma:=\inf\{\mu([y_i,z_i))\,,\,i=1,\ldots M_\gamma\}>0$. Furthermore,
if $\varepsilon$ is small enough,
$$\supp\mu_\varepsilon\subset \bigcup_{i=1}^{M_\gamma}[y_i,z_i)$$
and
$$\card\{j:x^\varepsilon_j\in[y_i,z_i)\}>{a_\gamma\over 2\varepsilon}\geq 1,$$
for every $i=1,\ldots, M_\gamma$. {}From now on, we assume that
$\varepsilon\in \cE$
is small enough so that the last two conditions hold, and we set
$$n^\varepsilon_i=\inf\{j:x^\varepsilon_j\in[y_i,z_i)\},\qquad i=1,\ldots
M_\gamma. $$
Denote by $\wt\tau^\varepsilon_k$ the $k$-th return of
$\wt\beta^\varepsilon$ to the origin. We also set
$\sigma^\varepsilon_i:=\wt\tau^\varepsilon_{n^\varepsilon_i}$
and $\sigma^\varepsilon_{M_\gamma+1}:=\wt\tau^\varepsilon_{N_\varepsilon}=
\tau^\varepsilon$.

Note that each of the variables $\sigma^\varepsilon_{i+1}-\sigma^\varepsilon_i$
is bounded below in distribution by $\wt
\tau^\varepsilon_{[a_\gamma/2\varepsilon]}$,
and recall that for every $c>0$, $\wt\tau^\varepsilon_{[c/\varepsilon]}$
converges in distribution to $\tau_c$. Since $\tau_c>0$ a.s., we may choose
$\delta\in(0,\kappa/2)$ so small that, for $\varepsilon$ small,
$$P[\sigma^\varepsilon_{i+1}-\sigma^\varepsilon_i> 2\delta
\hbox{ for every }i\in\{1,\ldots,M_\gamma\}]>1-{\eta\over 5}.
\eqno{(4.11)}.$$

Write $E'_\varepsilon$
for the intersection of the set $E_\varepsilon$
with the event considered in (4.11). Notice that on $E'_\varepsilon$
we can choose a finite sequence $0=s^\varepsilon_0<s^\varepsilon_1<\cdots
<s^\varepsilon_{K_\varepsilon}=\tau^\varepsilon$
in such a way that $\delta\leq s^\varepsilon_j-s^\varepsilon_{j-1}\leq
2\delta<\kappa$,
for every $j\in\{1,\ldots,K_\varepsilon\}$, and each interval
$[s^\varepsilon_{j-1},s^\varepsilon_j)$ is contained in
exactly one interval
$[\sigma^\varepsilon_{k-1},\sigma^\varepsilon_k)$.

We use the sequence $(s^\varepsilon_i)$ to get an upper bound
on $\theta(\varepsilon,\delta)$ on the event $E'_\varepsilon$. First
observe that
for $j\in\{1,\ldots,K_\varepsilon\}$,
$$\sup_{s,s'\in[s^\varepsilon_{j-1},s^\varepsilon_j)} d(\wt W^\varepsilon_s,\wt
W^\varepsilon_{s'})
\leq \sup_{s,s'\in[s^\varepsilon_{j-1},s^\varepsilon_j)} |\wt
\beta^\varepsilon_s-
\wt \beta^\varepsilon_{s'}|
+\sup_{s,s'\in[s^\varepsilon_{j-1},s^\varepsilon_j)} \sup_{t\geq 0}|\wt
W^\varepsilon_s(t)-\wt W^\varepsilon_{s'}(t)|.$$
The first term on the right hand side is bounded above by $\rho\leq \eta/5$ by
the definition of
$E_\varepsilon$ (cf (4.10)) and the property
$s^\varepsilon_j-s^\varepsilon_{j-1}
<\kappa$. To bound the second term, let
$s,s'\in[s^\varepsilon_{j-1},s^\varepsilon_j)$
and consider first the case when
$$m^\varepsilon(s,s'):=\inf_{r\in[s,s']}\wt\beta^\varepsilon_r>0.$$
Then $\wt W^\varepsilon_s(t)=\wt W^\varepsilon_{s'}(t)$
for every $t\in[0,m^\varepsilon(s,s')]$, and thus
$$\eqalign{&\sup_{t\geq 0}\,|\wt
W^\varepsilon_s(t)-\wt W^\varepsilon_{s'}(t)|\cr
&\quad\leq \sup_{m^\varepsilon(s,s')\leq t\leq \wt\beta^\varepsilon_s}|\wt
W^\varepsilon_s(t)-
\wt W^\varepsilon_s(m^\varepsilon(s,s'))|
+\sup_{m^\varepsilon(s,s')\leq t\leq \wt\beta^\varepsilon_{s'}}|\wt
W^\varepsilon_{s'}(t)-
\wt W^\varepsilon_{s'}(m^\varepsilon(s,s'))|\cr
&\quad\leq {2\eta\over 5}}$$
again by the definition of $E_\varepsilon$ (cf (4.9) and (4.10)). The case
$m^\varepsilon(s,s')=0$ is analogous, but we now get the additional
term $|\wt W^\varepsilon_s(0)-\wt W^\varepsilon_{s'}(0)|$. However, by
construction,
$s$ and $s'$ belong to the same interval
$[\sigma^\varepsilon_{k-1},\sigma^\varepsilon_k)$
and thus $\wt W^\varepsilon_s(0)$ and $\wt W^\varepsilon_{s'}(0)$
belong to the same $[y_k,z_k)$, which implies that $|\wt W^\varepsilon_s(0)-\wt
W^\varepsilon_{s'}(0)|\leq \gamma=\eta/5$. Finally, for every
$j\in\{1,\ldots,K_\varepsilon\}$, we get the bound
$$\sup_{s,s'\in[s^\varepsilon_{j-1},s^\varepsilon_j)} d(\wt W^\varepsilon_s,\wt
W^\varepsilon_{s'})\leq {4\eta\over 5}<\eta$$
on $E'_\varepsilon$. It follows that, for $\varepsilon$ small,
$$P[\theta(\varepsilon,\delta)\geq \eta]\leq P[(E'_\varepsilon)^c]\leq
{3\eta\over
5}<\eta.$$
This completes the proof of (4.8) and of the tightness of the sequence $\wt
W^\varepsilon$.

The remaining assertions of Proposition 4.2 are easy. (i) is clear since
$\wt \beta$ must be the weak limit of $\wt \beta^\varepsilon$. (ii)
follows from the analogous property
for $\wt W^\varepsilon$, and a similar argument applies to (iii).
\hfill$\square$

\medskip
\noi{\bf 4.3 \ Tightness of the reflected historical processes}
\medskip
Recall that the historical process for the $\varepsilon$-reflected system
is the process with values in $M_f({\cal W})$ defined by
$$\wt Y^\varepsilon_t=\int_0^{\tau_\varepsilon} d\wt
L^{\varepsilon,t}_s\,\delta_{
\wt W^\varepsilon_s}.$$
It is easy
to verify that $\wt Y^\varepsilon$ has right-continuous paths
with left limits. The following theorem is a slightly more precise
version of Theorem 1.1.

\smallskip
\noi{\bf Theorem 4.3}. {\it The sequence of the laws
$\wt \cL_Y^\eps$ of $\wt Y^\varepsilon$,
$\varepsilon\in\cE$, is tight in the space of probability
measures on $\D([0,\infty),M_f({\cal W}))$ and any limit law is
supported on $\C([0,\infty),M_f({\cal W}))$.
Suppose that $\wt \cL_Y$
is the limit of a subsequence of $\wt \cL_Y^\varepsilon$.
By passing to a further subsequence of $\eps$'s, if necessary, we may
assume that
the laws $\wt \cL_W^\eps$ of $\wt W^\varepsilon$
converge to a law $\wt \cL_W$.
Then one can construct on some probability space
processes $\wt Y$ and $\wt W$ with distributions
$\wt \cL_Y$ and $\wt \cL_W$, resp., related by
$$\wt Y_t=\int_0^{\tilde \tau} d\wt L^t_s\,\delta_{\wt W_s}\,,$$
where $(\wt L^t_s,t\geq 0,s\geq 0)$ denote the local times
of the process $\wt \beta_s:=\zeta_{\wt W_s}$, and
$\wt \tau=\inf\{s\geq 0:\wt L^0_s=a\}$.}
\medskip
\noi{\bf Proof}. By Proposition 4.2, the laws of
$\wt W^\varepsilon$, $\varepsilon\in \cE$ are tight. Hence, from any
subsequence
of $\cE$, we can extract a further subsequence $\cE_0$ along which
$\wt W^\varepsilon$ converges in distribution. We can in fact obtain more.
For every $\varepsilon>0$ and $t\geq 0$, denote by
$\Gamma^\varepsilon_t$, $\wt \Gamma^\varepsilon_t$ the random measures on
$\R_+$
defined by
$$\langle \Gamma^\varepsilon_t,\varphi\rangle=\int_0^{\tau\varepsilon}
dL^{\varepsilon,t}
_s\,\varphi(s)
\ ,\qquad \langle
\wt\Gamma^\varepsilon_t,\varphi\rangle=\int_0^{\tilde\tau\varepsilon}
d\wt L^{\varepsilon,t} _s\,\varphi(s).$$
(We have $\tau^\varepsilon=\wt \tau^\varepsilon$ but we prefer
to keep a different notation here.) Also define $\Gamma_t$ by:
$$\langle \Gamma_t,\varphi\rangle=\int_0^{\tau} dL^{t}
_s\,\varphi(s).$$
As a consequence
of Lemma 2.1, we know that
$$(\beta^\varepsilon_{\cdot\wedge \tau^\varepsilon},\Gamma^\varepsilon)
\build\hbox to 10mm{\rightarrowfill}
_{\varepsilon\to 0}^{}(\beta_{\cdot\wedge\tau},\Gamma)$$
uniformly on $[0,\infty)^2$, a.s. If we
replace the pair $(\beta^\varepsilon_{\cdot\wedge
\tau^\varepsilon},\Gamma^\varepsilon)$
by $(\wt \beta^\varepsilon_{\cdot\wedge
\tilde\tau^\varepsilon},\wt\Gamma^\varepsilon)$
this convergence still holds
in distribution in $\C(\R_+,\R)\times \D(\R_+,M_f(\R_+))$.
{}From this observation and
standard arguments, we have the joint convergence
$$(\wt W^\varepsilon,\wt\beta^\varepsilon,\wt\Gamma^\varepsilon)
\build\hbox to 10mm{\rightarrowfill}
_{\varepsilon\to 0,\varepsilon\in\cE_0}^{\rm (d)} (\wt W,\wt \beta,\wt \Gamma)
\eqno{(4.12)}$$
where
$$\langle \wt\Gamma_t,\varphi\rangle=\int_0^{\tilde \tau} d\wt L^{t}
_s\,\varphi(s),$$
with the notation introduced in the theorem.

By the Skorohod representation theorem, we can replace for every
$\varepsilon\in\cE_0$ the triplet $(\wt
W^\varepsilon,\wt\beta^\varepsilon,\wt\Gamma^\varepsilon)$
by a new triplet having the same distribution, in such a way that
the convergence (4.12) now holds almost surely. Without
risk of confusion, we keep the same notation for the new triplets.
We claim that we have then
$$\wt Y^\varepsilon_t=\int \wt\Gamma^\varepsilon_t(ds)\,\delta_{\wt
W^\varepsilon_s}
\build\hbox to 12mm{\rightarrowfill}
_{\varepsilon\to 0,\varepsilon\in\cE_0}^{}
\int \wt\Gamma_t(ds)\,\delta_{\wt W_s}=\wt Y_t\eqno{(4.13)}$$
uniformly on compact subsets of $\R_+$, a.s. Clearly Theorem 4.3 follows from
(4.13) and the fact that the limiting process
$\wt Y$ that appears in (4.13) is continuous. Both (4.13)
and the latter fact are
immediate consequences of the convergence
(4.12) (now assumed to hold
a.s.) and the following ``elementary'' lemma, whose
proof is left to the reader.
\medskip
\noi{\bf Lemma 4.4}. {\it Let $(\gamma^n,n\in\N)$
be a sequence in $\D(\R_+,M_f(\R_+))$. Assume that $\gamma^n_t$
converges as $n\to \infty$ to $\gamma_t$, uniformly
on every compact of $\R_+$, that $t\to\gamma_t$ is continuous and
that the measure $\gamma_t$ is diffuse, for every $t\in\R_+$.
Let $E$ be a Polish space and let $(f_n,n\in \N)$ be
a sequence in $\D(\R_+,E)$ that converges to $f$ in $\D(\R_+,E)$.
For every integer $n\in\N$ and every $t\in\R_+$,
let $\nu^n_t\in M_f(E)$ be defined by
$$\nu^n_t=\int \gamma^n_t(ds)\,\delta_{f_n(s)}.$$
Then $\nu^n_t$ converges as $n\to\infty$, uniformly
on compact subsets of $\R_+$, to the measure $\nu_t$
defined by
$$\nu_t=\int \gamma_t(ds)\,\delta_{f(s)}.$$
Furthermore, the mapping $t\to\nu_t$ is continuous.}
\hfill$\square$
\medskip
\noi{\bf Remark}. We do not know whether the limit law of
the sequence $\widetilde{\cal L}^\varepsilon_Y$ in Theorem 4.3
is unique. A positive answer would give the convergence
in distribution of the processes $\widetilde Y^\varepsilon$.
We can also formulate the problem in terms of the reflected
snake. Is there a unique (in law) process $\widetilde W$
satisfying properties (i) -- (iii) of Proposition 4.2 and such that
$$t\longrightarrow \int_0^{\tilde \tau} d\widetilde
L^t_s\,\delta_{\widetilde W_s(t)}$$
is a super-Brownian motion started at $\mu$ ?
\bigskip
\noi{\bf 5. Path properties of the reflected historical process}
\medskip
\noi{\bf 5.1 Preliminaries}
\medskip
Throughout this section, we consider a process $\wt Y$ which is a
weak limit of the processes $\wt Y^\varepsilon$ as $\varepsilon\to 0$.
According to Theorem 4.3, we may and will assume that
$\wt Y$ is constructed together with the reflected Brownian snake $\wt W$,
in such a way that, for every $t\geq 0$,
$$\wt Y_t=\int_0^{\tilde \tau} d\wt L^t_s\,\delta_{\wt W_s}$$
where $(\wt L^t_s,t\geq 0,s\geq 0)$ denote the local times
of the process $\wt \beta_s:=\zeta_{\wt W_s}$, which is (twice)
a reflected Brownian motion stopped at time
$\wt \tau=\inf\{s\geq 0:\wt L^0_s=a\}$.

The process
$$X_t=\int_0^{\tilde \tau} d\wt L^t_s\,\delta_{\wt W_s(t)}$$
is the weak limit of the processes $X^\varepsilon=\wt X^\varepsilon$
and therefore must be a super-Brownian motion started at $\mu$.

Let us recall the two key properties of the reflected snake $\wt W$
(cf Proposition 4.2):
\smallskip
$\bullet$ {\it Monotonicity property}: Almost surely for every $s\leq s'$
we have $\wt
W_s(t)\leq
\wt W_{s'}(t)$ for every $t\in[0,\wt\beta_s\wedge\wt\beta_{s'}]$.
\smallskip
$\bullet$ {\it Snake property}: The set of discontinuities of the mapping
$s\to \wt W_s$
is contained in the zero set of $\wt \beta$. Furthermore, if
$s<s'$ belong to the same connected component of the complement of
the zero set, we have
$$\wt W_s(t)=\wt W_{s'}(t)\qquad\hbox{for every
}t\in[0,\inf_{r\in[s,s']}\wt\beta_r]\,.$$

In order to state a useful preliminary result,
we introduce some notation. Let us fix $t>0$, and denote by
$(a^t_i,b^t_i)$, $i\in I_t$ the excursion intervals of
$\wt \beta$ above level $t$ (equivalently, these are the
connected components of the open set $\{s\geq 0:\wt \beta_s>t\}$).
Note that the index set $I_t$ may be empty. For each $i\in I_t$,
denote by $e^t_i$ the corresponding excursion
$$e^t_i(s)=\wt \beta_{(a^t_i+s)\wedge b^t_i}-t\,,\qquad s\geq 0.$$
By the snake property of $\wt W$, we have
$$\wt W_s(t)=\wt W_{a^t_i}(t)=:z^t_i\,,\quad \forall s\in[a^t_i,b^t_i].$$

We denote by $n(de)$ the
It\^o measure of positive Brownian excursions.
We normalize the measure $n(de)$ by declaring that
the Poisson point process of excursions from $0$,
i.e., the family of points $(L^0_{a^0_i}, e^0_i)$,
has intensity $ds\,n(de)$.
\medskip
\noi{\bf Proposition 5.1}. {\it Conditionally on $X_t$, the
point measure
$$\sum_{i\in I_t} \delta_{(z^t_i,e^t_i)}$$
is Poisson with intensity $X_t(dz)\,n(de)$. Consequently,
for every Borel subset $A$ of $\R$, the process
$$r \to Z^{t,A}_r=\int \wt Y_{t+r}(dw)\,{\bf 1}_A(w(t))$$
is a Feller diffusion started at $X_t(A)$.}
\medskip
We recall that the Feller diffusion is a
diffusion process $Z$ on $\R_+$ whose transition
kernels are characterized by the Laplace transform:
$E[\exp(-\lambda Z_t)| Z_0=z]=\exp(-z\,u_t(\lambda))$
where
$$u_t(\lambda)={\lambda \over 1+{1\over
2}\lambda t}.$$ The total mass process $\langle X_t,1\rangle=\wt
L^t_{\tilde \tau}$ is a Feller
diffusion started at $a$.

\medskip
\noi{\bf Proof}. We denote by $\tau^{(t)}_r$ the right-continuous inverse
of the function $r\to\wt L^t_r$. Note that $\tau^{(t)}_r<\infty$
iff $r<\wt L^t_{\tilde \tau}=\langle X_t,1\rangle$. We can rewrite the
definition
of $X_t$ as
$$\langle X_t,\varphi\rangle=\int_0^{\wt L^t_{\tilde \tau}} dr\,
\varphi(\wt W_{\tilde\tau^{(t)}_r}(t)).\eqno{(5.1)}$$
We also set for every $r\geq 0$,
$$A^{(t)}_r=\int_0^r du\,{\bf 1}_{\{\wt \beta_u>t\}}$$
and we let $\gamma^{(t)}_r$ be the right-continuous inverse
of the function $r\to A^{(t)}_r$. Finally we set
$\wt \beta^{(t)}_r=\wt \beta_{\gamma^{(t)}_r}-t$, for
every $r\in[0,A^{(t)}_{\tilde\tau})$.

We then claim that, conditionally on $\{\wt L^t_{\tilde \tau}=x\}$, the
process $(\wt \beta^{(t)}_r,0\leq
r<A^{(t)}_{\tilde\tau})$ is a reflected Brownian motion started at $0$ and
killed at
the first hitting time of $x$ by its local time at level $0$, and is
independent
of the process $(\wt W_{\tilde\tau^{(t)}_r}(t),0\leq r<\wt L^t_{\tilde \tau})$.
Except for the independence statement, this is a familiar property of linear
Brownian motion: See e.g. Section VI.2 of [RY]. To get the independence
property,
observe that the analogue of the process $\wt \beta^{(t)}$ for
the $\varepsilon$-reflected system codes (in the sense of Section 2) the
genealogy of the
descendants of particles at time $t$. On the other hand,
if $\tau^{\varepsilon,(t)}$ denotes the right-continuous inverse
of $\wt L^{\varepsilon,t}$, the process $(\wt
W^\varepsilon_{\tau^{\varepsilon,(t)}_r}(t),r\geq 0)$
just enumerates in increasing order the positions of the particles alive at
$t$. The
required independence is thus clear at the discrete level of the
$\varepsilon$-reflected system, and it is preserved under the
passage to the limit (4.12).

To complete the proof, write $\ell^{(t)}_i$ for the
local time at $0$ of $\wt\beta^{(t)}$ at the beginning, or the end, of
excursion
$e^t_i$. Note that $\tau^{(t)}_{\ell^{(t)}_i}=b^t_i$ and thus
$$z^{(t)}_i=\wt W_{\tau^{(t)}_{\ell^{(t)}_i}}(t).\eqno{(5.2)}$$
The point measure $\sum \delta_{(\ell^{(t)}_i,e^t_i)}$
is the excursion process of the process $\wt \beta^{(t)}$.
Hence,
conditionally on $\{\wt L^t_{\tilde \tau}=x\}$, this point measure is Poisson
with intensity $1_{[0,x)}(\ell)d\ell\,n(de)$
and is independent of $(\wt W_{\tilde\tau^{(t)}_r}(t),0\leq r<
\wt L^t_{\tilde \tau})$. The first part of the
lemma then follows from this property, (5.2) and (5.1) (which just says
that $X_t$ is
the image of the measure $1_{[0,\wt L^{(t)}_{\tilde \tau})}(\ell)d\ell$
under the mapping $\ell\to \wt W_{\tau^{(t)}_{\ell}}(t)$).

To get the second assertion of the lemma, note that by the definition of
$\wt Y_{t+r}$,
$$Z^{t,A}_r=\sum_{i\in I_t} {\bf 1}_{\{z^{(t)}_i\in A\}}\,\ell^r(e^{(t)}_i),$$
where $\ell^r(e^{(t)}_i)$ denotes the total local time of excursion $e^{(t)}_i$
at level $r$. By the first part of the proposition, conditionally on $X_t$,
the random measure
$$\sum_{i\in I_t} {\bf 1}_{\{z^{(t)}_i\in A\}}\,\delta_{e^{(t)}_i}$$
is Poisson with intensity $X_t(A)\,n(de)$.
Hence, conditionally
on $\{X_t(A)=x\}$, the process $(Z^{t,A}_r,r\geq 0)$
has the same law as $(L^r_{\tau_x},r\geq 0)$, and the desired
result follows from the celebrated Ray-Knight theorem
on Brownian local time. \hfill$\square$

\smallskip
\noi{\bf Remark}. We could easily sharpen the statement of Proposition 5.1 by
conditioning on $\wt Y_t$, or even on $(\wt Y_u,u\leq t)$ rather than on
$X_t$.
We will not need these refinements.

\medskip
\noi {\bf 5.2 \ A priori estimates}
\medskip
By [KS] or [R], we know that,
almost surely for every $t>0$, the measure $X_t$ has a continuous density
$x_t(y)$ with
respect to Lebesgue measure on $\R$, and the family
$(x_t(y),t>0,y\in\R)$ is jointly continuous.
Some of our results will be proved under the following
additional assumption:

\smallskip
\noi{\bf Assumption (H)}. {\it The measure $\mu$ has a continuous
density $x_0(y)$ with respect to Lebesgue measure.}
\smallskip
\noi Under (H), the family $(x_t(y),t\geq 0,y\in\R)$ is jointly continuous
(see Theorem 8.3.2 in [Da]).

In order to simplify the statements of the results in this
subsection we introduce a constant $\alpha$. All the results
hold for $\alpha =0$, assuming (H). Without this assumption,
the results hold for any fixed strictly positive $\alpha$.

\smallskip
For every $t\geq 0$, $r>0$ and $z\in \R$, we set
$$\psi_{t,t+r}(z)=\sup\{\wt W_s(t+r):\wt\beta_s\geq t+r \hbox{ and }\wt
W_s(t)<z \},$$
with the usual convention $\sup\emptyset=-\infty$.
We also consider the symmetric quantity:
$$\wh\psi_{t,t+r}(z)=\inf\{\wt W_s(t+r):\wt \beta_s\geq t+r \hbox{ and }
\wt W_s(t)>z
\},$$

\smallskip
\noi{\bf Proposition 5.2}. {\it Let $\eta\in(0,{1\over 2})$ and $c>0$.
Then, almost surely,
one can choose $\delta_0>0$ small enough so that, for every
$\delta\in(0,\delta_0)$,
$t\geq \alpha$ and $z\in\R$, the condition $x_t(z)\geq c$ implies
$$\psi_{t,t+\delta}(z)\geq z-\delta^{{1\over 2}-\eta}.$$}

\noi{\bf Proof}. For every $t\geq 0$ and $z\in\R$ set
$$\gamma^{t,z}=\inf\{s\geq 0:\wt\beta_s\geq t \hbox{ and } \wt W_s(t)\geq
z\},$$
with the convention $\inf\emptyset=\wt\tau$. Using the formula for $X_t$
in terms of $\wt W$,
and then the monotonicity property, we get
$$X_t((-\infty,z])=\int_0^{\tilde\tau} d\wt L^t_s\,{\bf 1}_{\{\wt W_s(t)<z\}}
=\int_0^{\tilde\tau} d\wt L^t_s\,{\bf 1}_{\{s<\gamma^{t,z}\}}=\wt
L^t_{\gamma^{t,z}}.$$
On the other hand, if $s<\gamma^{t,z}$ and $\wt \beta_{s}\geq t+\delta$,
we have $\wt W_s(t)<z$ and $\wt W_s(t+\delta)\leq \psi_{t,t+\delta}(z)$.
Therefore,
$$X_{t+\delta}((-\infty,\psi_{t,t+\delta}(z)])
=\int_0^{\tilde\tau} d\wt L^{t+\delta}_s\,
{\bf 1}_{\{\wt W_s(t+\delta)\leq \psi_{t,t+\delta}(z)\}}\geq
\wt L^{t+\delta}_{\gamma^{t,z}}.$$
Thanks to the H\"older continuity of Brownian local time in the time variable,
we can choose $\delta_1>0$ so small that, for every $\delta\in(0,\delta_1]$,
$t\geq 0$ and $z\in\R$,
$$\wt L^{t+\delta}_{\gamma^{t,z}}\geq \wt L^{t}_{\gamma^{t,z}}-
\delta^{{1\over 2}-\eta}.$$
By combining all these facts we obtain for every
$\delta\in(0,\delta_1]$,
$t\geq 0$ and $z\in\R$,
$$X_{t+\delta}((-\infty,\psi_{t,t+\delta}(z)])
\geq X_t((-\infty,z])-\delta^{{1\over 2}-\eta}.\eqno{(5.3)}$$
Note that the set $\{(t,y):x_t(y)>0\}$ is contained in the graph of $X$
and is thus relatively compact.
By uniform continuity, we
can choose $\delta_2>0$ small enough so that, for every $t\geq \alpha$
and $z\in\R$, the condition $x_t(z)\geq c$ implies that
$x_{t+\delta}(y)>{c\over 2}$ for all $\delta\in[0,\delta_2]$
and $y\in[z-\delta_2,z+\delta_2]$. In particular, if $0<r<\delta_2$ and
$\delta\in[0,\delta_2]$,
$$X_{t+\delta}((-\infty,z-r])<X_{t+\delta}((-\infty,z])-{c\over 2}\,r.
\eqno{(5.4)}$$
The proof of the following simple estimate for super-Brownian motion
is postponed to the appendix.
\smallskip
\noi {\bf Lemma 5.3}. {\it Almost surely there exists $\delta_3>0$ such that,
for every $t\geq \alpha$, $z\in\R$ and $\delta\in(0,\delta_3)$,
$$|X_{t+\delta}((-\infty,z])- X_t((-\infty,z])|\leq\delta^{{1\over 2}-\eta}.
\eqno{(5.5)}$$}

To complete the proof of Proposition 5.2, choose $\delta_0\in(0,\delta_1\wedge
\delta_2\wedge \delta_3)$ and also such that ${4\over c}\delta_0^{{1\over
2}-\eta}<\delta_2$.
Then, if $t\geq \alpha$ and $z\in\R$ are such that $x_t(z)\geq c$,
(5.3) and (5.5) give for $\delta\in(0,\delta_0)$,
$$X_{t+\delta}((-\infty,\psi_{t,t+\delta}(z)])
\geq X_{t+\delta}((-\infty,z]) - 2\,\delta^{{1\over 2}-\eta}.$$
Using (5.4) with $r={4\over c}\,\delta^{{1\over 2}-\eta}$, we get
$$X_{t+\delta}((-\infty,\psi_{t,t+\delta}(z)])
>X_{t+\delta}((-\infty,z-{4\over c}\,\delta^{{1\over 2}-\eta}]),$$
which implies
$$\psi_{t,t+\delta}(z)\geq z-{4\over c}\,\delta^{{1\over 2}-\eta}.$$
By replacing $\eta$ with $\eta'\in(0,\eta)$ we can get rid of the
factor ${4\over c}$.\hfill$\square$

\smallskip
We can immediately use Proposition 5.2 to derive some useful
results on continuity properties of the paths
$\wt W_s$. Note that, if $s\in(0,\wt\tau)$
is such that $\wt \beta_s\geq t+r$ and $\wt W_s(t)\geq z$, the monotonicity
property of the reflected snake implies that $\wt W_s(t+r)\geq
\psi_{t,t+r}(z)$.
Using Proposition 5.2 and the symmetric result for
$\wh \psi_{t,t+r}(z)$, we get the following corollary.
Recall that we take $\alpha=0$
if {\rm (H)} is assumed to hold and $\alpha>0$
otherwise.
\smallskip

\noi{\bf Corollary 5.4}. {\it
Let $\eta\in(0,{1\over 2})$ and $c>0$. Then almost
surely we can choose $\delta_0$ small enough so that, for every $t\geq \alpha$
and every $s\in(0,\wt\tau)$ such that $\wt\beta_s>t$ and $x_t(\wt
W_s(t))\geq c$,
we have for every $r\in[t,(t+\delta_0)\wedge \wt\beta_s]$,
$$|\wt W_s(r)-\wt W_s(t)|\leq (r-t)^{{1\over 2}-\eta}.$$}

\noindent{\bf 5.3 \ The key technical lemma}
\medskip
Our aim is to refine the a priori estimates that were derived in
the previous subsection. To this end, we will need a crucial
technical lemma (Lemma 5.7 below), whose proof requires
coming back to the approximating branching particle systems.
Recall the notation $(W^\varepsilon,\beta^\varepsilon,Y^\varepsilon)$ of
the previous sections. A much simplified version of the arguments
of Section 4 yields the convergence in distribution
$$(W^\varepsilon,Y^\varepsilon)
\build\hbox to 10mm{\rightarrowfill}
_{\varepsilon\to 0}^{\rm (d)} (W,Y),
$$
where $W$ is a minor modification of the Brownian snake of [L2]
(to be precise, $W$ is obtained by concatenating a Poisson point process
of Brownian snake excursions with intensity $\int \mu(dy)\,\N_y\,$,
in the notation of [L2])
and $Y$ is the historical
super-Brownian motion
connected to $W$ via the formula
$$Y_t=\int_0^\tau dL^t_s\,\delta_{W_s},$$
where $(L^t_s,t\geq 0,s\geq 0)$ are the local times
of the lifetime process $\beta_s=\zeta_{W_s}$, which is
a reflected Brownian motion stopped at time
$\tau=\inf\{s\geq 0:L^0_s=a\}$. (Our notation is slightly
inconsistent with the previous sections, where $\beta$
was not stopped, but this should cause no confusion.)

On the other hand (cf the proof of Theorem 4.3), we may
and will assume that there is a sequence ${\cal E}_0$ of values
of $\varepsilon$
such that
$$(\wt W^\varepsilon,\wt Y^\varepsilon)
\build\hbox to 10mm{\rightarrowfill}
_{\varepsilon\to 0,\varepsilon\in\cE_0}^{\rm (d)} (\wt W,\wt Y).
$$
By a compactness argument, and replacing the sequence $\cE_0$ by
a subsequence if necessary, we have also
$$(W^\varepsilon,Y^\varepsilon,\wt W^\varepsilon,\wt Y^\varepsilon)
\build\hbox to 10mm{\rightarrowfill}
_{\varepsilon\to 0,\varepsilon\in\cE_0}^{\rm (d)} (W,Y,\wt W,\wt Y).
$$
By the Skorohod representation theorem, we can
for every $\varepsilon\in\cE_0$ find a 4-tuple which has the same distribution
as $(W^\varepsilon,Y^\varepsilon,\wt W^\varepsilon,\wt Y^\varepsilon)$ (and
for which we
keep the same notation), in such a way that the previous convergence now
holds a.s.:
$$(W^\varepsilon,Y^\varepsilon,\wt W^\varepsilon,\wt Y^\varepsilon)
\build\hbox to 10mm{\rightarrowfill}
_{\varepsilon\to 0,\varepsilon\in\cE_0}^{\rm (a.s.)} (W,Y,\wt W,\wt Y).
\eqno{(5.6)}$$
{}From now on we will restrict our attention to
values of $\varepsilon$ in the sequence $\cE_0$
and assume that (5.6) holds.
{}From the equality $\wt X^\varepsilon= X^\varepsilon$, we also have
$$\int_0^\tau dL^t_s\,\delta_{W_s(t)}=\int_0^{\tilde\tau} d\wt
L^t_s\,\delta_{\wt W_s(t)}
=X_t,$$
and we see that $\tau$ coincides with $\wt \tau$.

We introduce the following more restrictive version of Assumption (H):
\smallskip
\noi{\bf Assumption (H')}. {\it The measure $\mu$ has a continuous
density $x_0(y)$, which is
H\"older continuous with exponent ${1\over 2}-\delta$, for
every $\delta>0$.}
\smallskip
\noindent As in the case of Assumption (H),
in order to be able to use a single statement for a result with or without
Assumption (H'),
we take $\alpha=0$ if (H') holds and otherwise we let
$\alpha$ be a fixed strictly positive constant. We also fix a constant
$c\in(0,1)$.
\smallskip
Let $\eta,\eta',\rho$ be three positive constants, with $0<\eta<\eta'<1/4$
and $\rho\in(0,{1\over 2})$. For every $\delta\in(0,1)$, we denote by
$E(\delta)$ the event on which the following three conditions hold.

\item{A.} For every $s\geq 0$, $t\in[0,\beta_s]$, and $r\in
[t,(t+\delta)\wedge \beta_s]$,
$$|W_s(r)-W_s(t)|\leq {1\over 2}\,(r-t)^{{1\over 2}-\eta}.$$

\item{B.} For every $t\geq \alpha$ and $s\geq 0$ such that $\wt \beta_s>t$
and $x_t(\wt W_s(t))\geq c$, we have for every $r\in[t,(t+\delta)\wedge
\wt\beta_s]$,
$$|\wt W_s(r)-\wt W_s(t)|\leq (r-t)^{{1\over 2}-\eta}.$$

\item{C.} For every $t\geq \alpha$, $z\in\R$, and $y\in[z-\delta^{{1\over
2}-\eta'},z+\delta^{{1\over 2}-\eta'}]$,
$$|x_t(z)-x_t(y)|\leq |z-y|^{{1\over 2}-\rho}.$$

\noindent Note that the sets $E(\delta)$ are decreasing in $\delta$.
We have $P[\bigcup_n E(2^{-n})]=1$. The fact that
properties A and B hold for $\delta$ small enough follows
from the H\"older continuity properties of the Brownian snake paths
(cf (2.7))  and Corollary 5.4 respectively.
For property C, see Theorem 8.3.2 in [Da] when $\alpha>0$. When
$\alpha=0$ (then (H') is in force), the desired H\"older
continuity of the densities is easily obtained from formula (8.3.5b) of [Da]
by using the techniques of [KS].

Throughout this subsection, we fix $\delta\in(0,1)$, $t\geq \alpha$
and $z\in\R$. We plan to improve the estimates obtained on
$\psi_{t,t+\delta}(z)$ in the previous subsection. We set
$$\gamma=\delta^{{1\over 2}-\eta'}$$
and we assume that $\delta$ has been chosen small enough so that
$\gamma>4\,\delta^{{1\over 2}-\eta}$. Then, for every $r\in[t,t+\delta]$,
we set
$$X^*_r=\int Y_r(dw)\,\bone_{\{w(t)\in(z-\gamma,z+\gamma)\}}\,\delta_{w(r)}.$$
The random measure $X^*_r$ corresponds, for the historical super-Brownian
motion $Y$,
to the contribution of those particles alive at time $r$ whose ancestor at time
$t$ lies in the interval $(z-\gamma,z+\gamma)$. Note that $X^*_t$ is
simply the restriction of $X_t$ to $(z-\gamma,z+\gamma)$.

Our goal is to compare $X^*_{t+\delta}((-\infty,\psi_{t,t+\delta}(z)])$
to $X^*_t((-\infty,z])$ in the same way as we compared
$X_{t+\delta}((-\infty,\psi_{t,t+\delta}(z)])$
to $X_t((-\infty,z])$ in (5.3) above. Unfortunately, the argument
has to be significantly more complicated.

We set for every $\varepsilon>0$,
$$\psi^\varepsilon_{t,t+\delta}(z)
=\sup\{\wt W^\varepsilon_s(t+\delta):\wt\beta^\varepsilon_s
\geq t+\delta\hbox{ and }\wt W^\varepsilon_s(t)<z\},$$
which represents for the $\varepsilon$-reflected system the
right-most position among those particles alive at time $t+\delta$
which are descendants of the particles located to the left of $z$ at time $t$.

\medskip
\noi{\bf Lemma 5.5}. {\it We have
$$\psi_{t,t+\delta}(z)=\lim_{\varepsilon\to 0} \psi^\varepsilon_{t,t+\delta}(z)
\qquad \hbox{a.s.}$$}

\noi{\bf Proof}. This is basically a consequence of the convergence of $\wt
W^\varepsilon$
towards $\wt W$, which entails the convergence of $\wt \beta^\varepsilon$
to $\wt \beta$.
We also use the fact that in the definition of $\psi_{t,t+\delta}(z)$, i.e.,
$$\psi_{t,t+\delta}(z)
=\sup\{\wt W_s(t+\delta):\wt\beta_s
\geq t+\delta\hbox{ and }\wt W_s(t)<z\},$$
we can replace the weak inequality $\wt\beta_s
\geq t+\delta$ by a strict one, and/or the strict inequality $\wt W_s(t)<z$
by a weak one. To justify this, note that:
\smallskip
\item{(a)} Almost surely, every $s$ such that $\wt \beta_s=t+\delta$ is the
limit
of a sequence $s_n$ such that $\wt \beta_{s_n}>t+\delta$ (simply because
$t+\delta$ cannot be a local maximum of $\wt \beta$).
\smallskip
\item{(b)} With probability 1,
there is no value of $s$ such that $\wt W_s(t)=z$ and $\wt
\beta_s\geq
t+\delta$ (this immediately follows from Proposition 5.1).
\smallskip
\noi We leave details to the reader. \hfill$\square$

\smallskip
We now introduce a different approximation of $\psi_{t,t+\delta}(z)$.
We consider in the (non-reflected) $\varepsilon$-system
those particles which are located at time $t$ in the interval
$(z-\gamma,z+\gamma)$, and the descendants of these particles after
time $t$. With this branching particle system (evolving
over the time interval $[t,\infty)$), we can associate a reflected
system
in the way explained in Subsection 3.1. We denote by
$\psi^{*,\varepsilon}_{t,t+\delta}(z)$
the position in this new reflected system of the right-most particle at time
$t+\delta$, among those particles which are descendants of the particles
located to the left of $z$ at time $t$.

For every $r>0$, we set $\un{x}(t,z,r)=\inf\{x_t(y):|y-z|\leq r\}$
and $\ov{x}(t,z,r)=\sup\{x_t(y):|y-z|\leq r\}$.

\smallskip
\noi{\bf Lemma 5.6}. {\it We have
$$P\Big[\Big(\limsup_{\varepsilon\downarrow 0}
\{\psi^{*,\varepsilon}_{t,t+\delta}(z)\not
=\psi^{\varepsilon}_{t,t+\delta}(z)\}\Big)
\cap E(\delta)\cap\{\un{x}(t,z,\delta^{1/2})\geq c\}\Big]\leq 2
\exp(-2c\delta^{-1/2}).$$}

\noi{\bf Proof}. We introduce the following events:
$$\Lambda^+=\{\exists s\geq 0:\wt\beta_s>t+\delta \hbox{ and } \wt
W_s(t)\in(z-\delta^{1/2},z)\},$$
and
$$\Lambda^-=\{\exists s\geq 0:\wt\beta_s>t+\delta \hbox{ and } \wt
W_s(t)\in(z,z+\delta^{1/2})\}.$$
We first verify that a.s.,
$$\Big(\Big(\limsup_{\varepsilon\downarrow 0}
\{\psi^{*,\varepsilon}_{t,t+\delta}(z)\not
=\psi^{\varepsilon}_{t,t+\delta}(z)\}\Big)
\cap E(\delta)\cap\{\un{x}(t,z,\delta^{1/2})\geq c\}\Big)
\subset (\Lambda^+\cap\Lambda^-)^c.\eqno{(5.7)}$$
Suppose that $\Lambda^+\cap\Lambda^-\cap
E(\delta)\cap\{\un{x}(t,z,\delta^{1/2})
\geq c\}$
holds. Then, there exists
$s_1\geq 0$ such that $\wt\beta_{s_1}>t+\delta$
and $W_{s_1}(t)\in(z-\delta^{1/2},z)$. {}From property B in the definition
of $E(\delta)$ we also have $|\wt W_{s_1}(r)-z|<2\delta^{{1\over 2}-\eta}$
for every $r\in[t,t+\delta]$. Similarly, there exists $s_2\geq 0$
such that $\wt\beta_{s_2}>t+\delta$,
$W_{s_2}(t)\in(z,z+\delta^{1/2})$ and $|\wt W_{s_2}(r)-z|<2\delta^{{1\over
2}-\eta}$
for every $r\in[t,t+\delta]$. By the convergence (5.6), the
same properties hold for $\varepsilon>0$ small enough, if we replace $\wt
W_{s_i}$ and
$\wt \beta_{s_i}$ by $\wt W^\varepsilon_{s_i}$ and
$\wt \beta^\varepsilon_{s_i}$ respectively.

On the other hand, by property A of the definition of $E(\delta)$
and the convergence (5.6), we have
also for $\varepsilon$ small enough, for every $s$ such that
$\beta^\varepsilon_s\geq t$ and every $r\in[t,(t+\delta)\wedge
\wt\beta^\varepsilon_s]$,
$$|W^\varepsilon_s(r)-W^\varepsilon_s(t)|\leq \delta^{{1\over 2}-\eta}.$$
In particular, if $s$ is such that $\beta^\varepsilon_s\geq t$ and
$|W^\varepsilon_s(t)-z|\geq \gamma\geq 4\,\delta^{{1\over 2}-\eta}$, we
have for
every $r\in[t,(t+\delta)\wedge \wt\beta^\varepsilon_s]$,
$$|W^\varepsilon_s(r)-z|>2\,\delta^{{1\over 2}-\eta}.$$

We have shown that, on the event $\Lambda^+\cap\Lambda^-\cap
E(\delta)\cap\{\un{x}(t,z,\delta^{1/2})\geq c\}$, provided that
$\varepsilon$ is small
enough:

\smallskip
\item{$\bullet$} There exist $s_1$ and $s_2$ such that
$\wt\beta_{s_1}>t+\delta$,
$\wt\beta_{s_2}>t+\delta$ and
$$\eqalign{&W^\varepsilon_{s_1}(t)\in(z-\delta^{1/2},z)\;,\qquad
|\wt W^\varepsilon_{s_1}(r)-z|<2\delta^{{1\over 2}-\eta},\ \forall
r\in[t,t+\delta]\cr
&W^\varepsilon_{s_2}(t)\in(z,z+\delta^{1/2})\;,\qquad
|\wt W^\varepsilon_{s_2}(r)-z|<2\delta^{{1\over 2}-\eta},\ \forall
r\in[t,t+\delta].}$$

\item{$\bullet$} For every $s\geq 0$ such that $\beta^\varepsilon_s\geq t$ and
$|W^\varepsilon_s(t)-z|\geq \gamma$,
$$|W^\varepsilon_s(r)-z|>2\,\delta^{{1\over 2}-\eta},\ \forall
r\in[t,(t+\delta)\wedge
\wt\beta^\varepsilon_s].$$

\noi These properties allow us to apply Lemma 3.1. In the context
of that lemma, the original system
is the $\varepsilon$-system
considered after time $t$, the new (restricted)
system
consists of the descendants of the particles which are located
at time $t$ in the interval $(z-\gamma,z+\gamma)$, and we take
$I=(z-2\delta^{{1\over
2}-\eta},z+2\delta^{{1\over 2}-\eta})$.
Lemma 3.1 and the previous properties imply that the
restrictions of the paths $\wt W^\varepsilon_{s_1}$ and
$\wt W^\varepsilon_{s_2}$ to $[t,t+\delta]$ still appear as restrictions
of reflected historical paths in the new system.
Note that in the definition of $\psi^{\varepsilon}_{t,t+\delta}(z)$,
respectively of $\psi^{*,\varepsilon}_{t,t+\delta}(z)$, we may
restrict our attention to those reflected historical paths between
times $t$ and $t+\delta$
in the original system, resp. in the new system,
whose value at time $t$ lies in the interval $[\wt W^\varepsilon_{s_1}(t),z)$
(this is so because of the monotonicity property of reflected
historical paths).
Any such path is bounded below and above by $\wt W^\varepsilon_{s_1}$
and $\wt W^\varepsilon_{s_2}$ respectively, on the time interval
$[t,t+\delta]$.
By Lemma 3.1 again, the class of paths that we need to consider is exactly the
same for both the original system and the new one.
This is enough to conclude that
$\psi^{*,\varepsilon}_{t,t+\delta}(z)=\psi^\varepsilon_{t,t+\delta}(z)$, and
we get our claim (5.7).

It follows from (5.7) that the probability considered in the lemma
is bounded above by
$$P[(\Lambda^+\cap\Lambda^-)^c\cap \{\un{x}(t,z,\delta^{1/2})\geq c\}].$$
By the construction of $\wt Y$, we have
$$\int_0^{\tilde\tau} d\wt L^{t+\delta}_s\,\bone_{\{\wt
W_s(t)\in(z-\delta^{1/2},z)\}}
=\int \wt Y_{t+\delta}(dw)\,\bone_{\{w(t)\in(z-\delta^{1/2},z)\}}.$$
Hence the event $\Lambda^+$
certainly holds if
$$\int \wt Y_{t+\delta}(dw)\,\bone_{\{w(t)\in(z-\delta^{1/2},z)\}}>0.$$
It follows that
$$P[(\Lambda^+)^c\cap \{\un{x}(t,z,\delta^{1/2})\geq c\}]
\leq P[\{\int \wt Y_{t+\delta}(dw)\,\bone_{\{w(t)\in(z-\delta^{1/2},z)\}}=0\}
\cap \{\un{x}(t,z,\delta^{1/2})\geq c\}],$$
and a similar bound holds if we replace $\Lambda^+$ by $\Lambda^-$. By
Proposition 5.1, the
last quantity is bounded above by the probability that a Feller diffusion
started
at $c\,\delta^{1/2}$ vanishes at time $\delta$. This probability is
equal to $\exp(-2c\delta^{-1/2})$, which completes the proof.
\hfill$\square$
\smallskip
We can now state the key lemma. We fix still another constant
$\eta''\in(\eta',1/4)$.
\medskip
\noi{\bf Lemma 5.7}. {\it There exist two positive constants $C$ and $\kappa$,
that depend only on $c,\eta,\eta',\eta''$ and $\rho$, such that
$$\eqalign{
&P\big[\{|X^*_{t+\delta}((-\infty,\psi_{t,t+\delta}(z)])
-X^*_t((-\infty,z])|>\delta^{{3\over
4}-\eta''}\}\cr
&\qquad\qquad\cap E(\delta)\cap \{\un{x}(t,z,\delta^{1/2})\geq c\}
\cap\{\ov{x}(t,z,\gamma)\leq c^{-1}\}\big]
\leq C\,\exp(-\delta^{-\kappa}).}$$}
\medskip
\noi{\bf Proof}. For every $r\in[t,t+\delta]$, set
$$X^{*,\varepsilon}_r=\int
Y^\varepsilon_{r}(dw)\,\bone_{\{w(t)\in(z-\gamma,z+\gamma)\}}
\delta_{w(r)},$$
which represents the contribution at time $r$ of the descendants
(in the non-reflected system) of particles which are located in
$(z-\gamma,z+\gamma)$
at time $t$. {}From the convergence of $Y^\varepsilon$ to $Y$, and the
fact that $\int Y_r(dw)\,\bone_{\{w(t)=z\pm \gamma\}}=0$, one
can easily show
that for every $r\in [t,t+\delta]$, the measures $X^{*,\varepsilon}_r$
converge weakly to $X^{*}_r$. In particular, a.s. for every $y\in\R$,
$$\lim_{\varepsilon\to 0} X^{*,\varepsilon}_{t+\delta}((-\infty,y])
=X^*_{t+\delta}((-\infty,y]).$$
{}From Lemma 5.5 and Lemma 5.6, we get that on the set $E(\delta)
\cap \{\un{x}(t,z,\delta^{1/2}\geq c\}$, we have the convergence
$$\lim_{\varepsilon\to 0} X^{*,\varepsilon}_{t+\delta}((-\infty,
\psi^{*,\varepsilon}_{t,t+\delta}(z)])
=X^*_{t+\delta}((-\infty,\psi_{t,t+\delta}(z)]),$$
except possibly on a set of measure at most $2\,\exp(-2c\delta^{-1/2})$.

However, by the definition of $\psi^{*,\varepsilon}_{t,t+\delta}(z)$,
and the monotonicity property of reflected systems, the
quantity $X^{*,\varepsilon}_{t+\delta}((-\infty,
\psi^{*,\varepsilon}_{t,t+\delta}(z)])$ is equal to $\varepsilon$
times the number of descendants at time $t+\delta$ of the particles
present at time $t$ in $(z-\gamma,z)$, for the $\varepsilon$-reflected system
constructed over the time interval $[t,\infty)$ from the particles
present at time $t$ in $(z-\gamma,z+\gamma)$. Since the law of the
branching evolution
is the same for the reflected system as for the original one, we see that
conditionally on $\{X^{*,\varepsilon}_t((-\infty,z])=\varepsilon\,k\}$,
the variable $X^{*,\varepsilon}_{t+\delta}((-\infty,
\psi^{*,\varepsilon}_{t,t+\delta}(z)])$ is distributed as
$\varepsilon\,Z^{(\varepsilon,k)}_\delta$, where $Z^{(\varepsilon,k)}$
denotes a Galton-Watson process
with critical binary branching at rate $\varepsilon^{-1}$
and initial value $k$. Recall that
$X^{*,\varepsilon}_t((-\infty,z])$
converges a.s. to $X^*_t((-\infty,z])$. By standard limit theorems for
Galton-Watson
processes,
$$(X^{*,\varepsilon}_t((-\infty,z]),X^{*,\varepsilon}_{t+\delta}((-\infty,
\psi^{*,\varepsilon}_{t,t+\delta}(z)]))
\build\hbox to 10mm{\rightarrowfill}
_{\varepsilon\to 0}^{\rm (d)} (X^{*}_t((-\infty,z]),U),$$
where conditionally on $X^{*}_t((-\infty,z])=u$, the variable $U$ is
distributed as the value at time $\delta$ of a Feller diffusion started at
$u$.

Note that $X^*_t((-\infty,z])=X_t((z-\gamma,z])$ and that on the set
$\{\ov{x}(t,z,\gamma)\leq c^{-1}\}$ we have $X^*_t((-\infty,z])\leq
c^{-1}\gamma
= c^{-1}\delta^{{1\over 2}-\eta'}$. Elementary estimates on the Feller
diffusion,
using only the form of the Laplace transform of the semigroup  (see the
appendix for very similar estimates)
show that
$$P\big[\{X^*_t((-\infty,z])\leq c^{-1}\delta^{{1\over 2}-\eta'}\}
\cap \{|U-X^*_t((-\infty,z])|\geq \delta^{{3\over 4}-\eta''}\}\big]
\leq C'\,\exp(-\delta^{-\kappa'}),$$
where the constants $C'$ and $\kappa'>0$ depend only on $c,\eta'$ and $\eta''$.

To complete the proof of the lemma, we write
$$\eqalign{
&P\big[\{|X^*_{t+\delta}((-\infty,\psi_{t,t+\delta}(z)])-X^*_t((-\infty,z])|
>\delta^{{3\over
4}-\eta''}\}\cr
&\hskip 5cm\cap E(\delta)\cap \{\un{x}(t,z,\delta^{1/2})\geq c\}
\cap\{\ov{x}(t,z,\gamma)\leq c^{-1}\}\big]\cr
&\leq 2\,\exp(-2c\delta^{-1/2})
+P\big[\big\{\liminf_{\varepsilon\to 0}|X^{*,\varepsilon}_{t+\delta}((-\infty,
\psi^{*,\varepsilon}_{t,t+\delta}(z)])
-X^*_{t}((-\infty,z])|>\delta^{{3\over
4}-\eta''}\big\}\cr
&\hskip 5cm\cap \{X^*_t((-\infty,z])\leq c^{-1}\delta^{{1\over
2}-\eta'}\}\big]\cr
&\leq 2\,\exp(-2c\delta^{-1/2})
+\liminf_{\varepsilon\to 0}P\big[\big\{|X^{*,\varepsilon}_{t+\delta}((-\infty,
\psi^{*,\varepsilon}_{t,t+\delta}(z)])
-X^*_{t}((-\infty,z])|>\delta^{{3\over
4}-\eta''}\big\}\cr
&\hskip 5cm\cap \{X^*_t((-\infty,z])\leq c^{-1}\delta^{{1\over
2}-\eta'}\}\big]\cr
&\leq 2\,\exp(-2c\delta^{-1/2})
+P\big[\{X^*_t((-\infty,z])\leq c^{-1}\delta^{{1\over 2}-\eta'}\}
\cap \{|U-X^*_t((-\infty,z])|\geq \delta^{{3\over 4}-\eta''}\}\big]\cr
&\leq 2\,\exp(-2c\delta^{-1/2}) + C'\,\exp(-\delta^{-\kappa'}).
}$$
\par\hfill$\square$
\medskip
\noi{\bf 5.4 \ The main result}
\medskip
We keep the notation introduced in the previous subsection. The reals
$t\geq \alpha$, $\delta\in (0,1)$ and $z\in\R$ are fixed for the moment.
\medskip
\noi{\bf Lemma 5.8}. {\it Assume that $\eta''>{3\over 2}\eta'+{1\over
2}\rho$. There exist
two constants
$\ov C$ and
$\ov\kappa>0$, that depend only on $c,\eta',\eta''$ and $\rho$, such that
$$P\big[\{X^*_{t+\delta}((-\infty,z])\geq X^*_t((-\infty,z])+\delta^{{3\over 4}
-\eta''}\}\cap E(\delta)\cap \{\ov{x}(t,z,\gamma)\leq c^{-1}\}\big]
\leq \ov C \exp(-\delta^{-\ov \kappa}).$$}

The proof of this lemma is an application of standard techniques in the
theory of super-Brownian motion. See the appendix for a detailed argument.

\medskip
\noi{\bf Proposition 5.9}. {\it Under the assumptions of Lemma 5.8,
there exist two constants $C_0$ and $\kappa_0>0$,
that depend only on $c,\eta,\eta',\eta''$ and $\rho$, such that
$$\eqalign{
&P\big[\{\psi_{t,t+\delta}(z)<z-{2\over c}\delta^{{3\over 4}-\eta''}\}\cap
E(\delta)\cr
&\qquad\qquad \cap \{c\leq \un{x}(t,z,\gamma)\leq \ov{x}(t,z,\gamma)\leq
c^{-1}\}
\cap\{\un{x}(t+\delta,z,\gamma)>c\}\big]\leq C_0\,\exp(-\delta^{-\kappa_0}).}
$$}
\medskip
\noi{\bf Proof}. Our argument is very similar to the proof of Proposition
5.2. We
will assume that the event
$E(\delta)\cap \{c\leq \un{x}(t,z,\gamma)\leq
\ov{x}(t,z,\gamma)\leq c^{-1}\}$ holds.
By Lemmas 5.7 and 5.8, we have on this set
$$X^*_{t+\delta}((-\infty,z])\leq X^*_t((-\infty,z])+\delta^{{3\over 4}-\eta''}
\leq X^*_{t+\delta}((-\infty,\psi_{t,t+\delta}(z)])+2\,\delta^{{3\over
4}-\eta''}
\eqno{(5.8)}$$
except possibly on a set of probability at most
$C\,\exp(-\delta^{-\kappa})+\ov{C}\,\exp(-\delta^{-\ov{\kappa}})$.

On the other hand, condition A in the definition of $E(\delta)$
(and the fact that $\gamma>4\,\delta^{{1\over 2}-\eta}$) ensures that
the measures $X^*_{t+\delta}$ and $X_{t+\delta}$ coincide over the
interval $(z-{\gamma\over 2},z+{\gamma\over 2})$. Hence, on the event
$\{\un{x}(t+\delta,z,\gamma)>c\}$, we get
$$X^*_{t+\delta}((-\infty,z])-2\,\delta^{{3\over 4}-\eta''}
>X^*_{t+\delta}((-\infty,z-{2\over c}\delta^{{3\over 4}-\eta''}]),$$
provided that $\delta$ is small enough so that ${2\over c}\delta^{{3\over
4}-\eta''}
\leq {\gamma\over 2}$. On the set where (5.8) holds, we get
$$X^*_{t+\delta}((-\infty,\psi_{t,t+\delta}(z)])
>X^*_{t+\delta}((-\infty,z-{2\over c}\delta^{{3\over 4}-\eta''}]),$$
and the desired result follows. \hfill$\square$

\medskip
We now come to the main result of this section, which is
a refinement of Corollary 5.4.
Recall our conventions concerning $\alpha$---this constant
is equal $0$ if (H') is assumed to hold and otherwise
$\alpha$ is a fixed strictly positive constant.
\medskip
\noi{\bf Theorem 5.10}. {\it
Let $\lambda>0$ and $c\in(0,1)$. Then a.s.
we can choose $\delta_0$ small enough so that, for every $t\geq \alpha$
and every $s\in(0,\tau)$ such that $\wt \beta_s>t$ and $x_t(\wt W_s(t))\geq c$,
we have for every $r\in[t,(t+\delta_0)\wedge \wt \beta_s]$,
$$|\wt W_s(r)-\wt W_s(t)|\leq (r-t)^{{3\over 4}-\lambda}.$$}

\noi{\bf Proof}. We can choose $\eta,\eta',\eta''$ with
$0<\eta<\eta'<\eta''<\lambda$
and $\rho\in(0,{1\over 2})$ such that the assumptions of Lemma 5.8 hold.
We then apply the estimate of Proposition 5.9 with $\delta=2^{-n}$ ($n$
large enough) to all reals $t\in[\alpha,n]$, $z\in[-n,n]$
of the form $t=k2^{-n}$, $z=j2^{-n}$. We have already observed that
$P[\bigcup_n E(2^{-n})]=1$. Furthermore, if we assume that $c\leq
x_t(z)\leq c^{-1}$
we will have $\un{x}(t,z,2^{-n({1\over
2}-\eta')})\geq c/2$,
$\ov{x}(t,z,2^{-n({1\over 2}-\eta')})\leq 2/c$,
and $\un{x}(t+2^{-n},z,2^{-n({1\over 2}-\eta')})\geq c/2$, for all $n$
sufficiently large
(depending on $\omega$ but not on $t$ and $z$). Then, by
combining the estimate of Proposition 5.9 with the
Borel-Cantelli lemma, we obtain
the following property: There exists an integer
$n_0(\omega)$ such that for every $n\geq n_0(\omega)$, for
every $t=k2^{-n}$, $z=j2^{-n}$ with $t\in[\alpha,n]$, $z\in[-n,n]$,
the condition $c\leq x_t(z)\leq c^{-1}$ implies
$$\psi_{t,t+2^{-n}}(z)\geq z-(2^{-n})^{{3\over 4}-\lambda}.$$
Since the densities $x_r(y)$ are bounded over $[\alpha,\infty)\times\R$,
a simple argument shows that we can drop the condition $x_t(z)\leq c^{-1}$
in the previous assertion.

Then, if $s\geq 0$ is such that $\wt \beta_s\geq t+2^{-n}$, where
$t$ is of the form $t=k2^{-n}$, we let
$z=j2^{-n}$ be such that $z<\wt W_s(t)\leq z+2^{-n}$. If $n$
is large enough (again independently of
the choice of $s$ and $t$), the condition $x_t(\wt W_s(t))\geq 2c$
will imply $x_t(z)>c$. Then, by the definition of $\psi_{t,t+\delta}(z)$
and the preceding estimate,
$$\wt W_s(t+2^{-n})\geq \psi_{t,t+2^{-n}}(z)\geq \wt W_s(t)-2^{-n}
-(2^{-n})^{{3\over 4}-\lambda}.$$
Thanks to this observation and a symmetry argument, we obtain that a.s.
for $n$ large enough, for every $t\geq \alpha$ of the form
$t=k2^{-n}$ and every $s\geq 0$ such that $\wt \beta_s\geq t+2^{-n}$
and $x_t(\wt W_s(t))\geq 2c$,
$$|\wt W_s(t+2^{-n})-\wt W_s(t)|\leq 2\,(2^{-n})^{{3\over 4}-\lambda}.$$
The statement of Theorem 5.10 now follows easily thanks to the usual
chaining argument. \par\hfill$\square$

\smallskip
Theorem 1.2 is an immediate consequence of Theorem 5.10. Note that,
by the representation formula for $\wt Y$ in terms of $\wt W$,
the set $\supp\wt Y_t$ is contained in
$\{\wt W_s;\wt \beta_s=t\}$, for every $t>0$, a.s. The comments
following the statement of Theorem 1.2 are justified
by Proposition 5.1.

\bigskip
\noi{\bf 6. Branching points}
\medskip
In this last section, we prove Theorem 1.3. As in Section 5, we assume
that the process $\wt Y$ is constructed together with the
reflected Brownian snake $\wt W$, in such a way that we have
the representation formula
$$\wt Y_t=\int_0^{\tilde \tau} d\wt L^t_s\,\delta_{\wt W_s}.$$
We need a preliminary lemma.
If $s_1<s_2$, we set $m(s_1,s_2)=\inf_{s\in[s_1,s_2]}\wt \beta_s$.
\medskip
\noi{\bf Lemma 6.1}. {\it Almost surely, for any $t>0$ and any $s_1<s_2$
such that
$\wt \beta_{s_1}=\wt \beta_{s_2}=t$ and $0<m(s_1,s_2)<t$, we have
$$x_{m(s_1,s_2)}(\wt W_{s_1}(m(s_1,s_2)))>0.$$}
\medskip
\noi{\bf Proof}. Let $\alpha>0$ and let $A\geq 1$ be an integer. Write
$E_A$ for the event
$E_A=\{{\cal G}\subset [0,A]\times [-A,A]\}$, where ${\cal G}$ is as above
the graph of $X$.
It is enough to prove that a.s. on $E_A$, the following holds:
\smallskip
\item{(P)} For any $t>\alpha$ and $s_1<s_2$ such that
$\wt \beta_{s_1}=\wt \beta_{s_2}=t$ and $0<m(s_1,s_2)<t-\alpha$, we have
$x_{m(s_1,s_2)}(\wt W_{s_1}(m(s_1,s_2)))>0$.
\smallskip
\noi We first introduce some notation. Let $e$ be an excursion, that is
a continuous function $e:\R_+\longrightarrow \R_+$ such that $e(s)>0$
iff $0<s<\sigma(e)$, for some $\sigma(e)>0$. Set
$$T_\alpha(e)=\inf\{s\geq 0:e(s)=\alpha\}$$
and, if $T_\alpha(e)<\infty$,
$$\eqalign{L_\alpha(e)&=\sup\{t\geq 0:e(t)=\alpha\},\cr
M_\alpha(e)&=\inf_{T_\alpha(e)\leq s\leq L_\alpha(e)}e(s)
}$$
By convention we take $M_\alpha(e)=0$ if $T_\alpha(e)=\infty$.

Let $r>0$. Recall the notation $I_r$ and
$e^r_i,z^r_i$, $i\in I_r$ introduced before Proposition 5.1, and for every
$c>0$ and
$\delta\in(0,\alpha)$, set
$$N^\delta_r(\alpha,c)=\sum_{i\in I_r} {\bf 1}_{\{x_r(z^r_i)\leq c\}}
\,{\bf 1}_{\{0<M_\alpha(e^r_i)\leq \delta\}}.$$
Proposition 5.1 allows us to conclude that,
$$\eqalign{E[N^\delta_r(\alpha,c)\,{\bf 1}_{E_A}]
&\leq E\Big[\sum_{i\in I_r} {\bf 1}_{\{|z^r_i|\leq A\}}\,{\bf
1}_{\{x_r(z^r_i)\leq c\}}
\,{\bf 1}_{\{0<M_\alpha(e^r_i)\leq \delta\}}\Big]\cr
&=E\Big[\int_{-A}^A dz\,x_t(z)\,{\bf 1}_{\{x_r(z)\leq
c\}}\,n(0<M_\alpha(e)\leq \delta)\Big]\cr
&\leq 2cA\,\delta\alpha^{-2},}$$
using the easy formula $n(0<M_\alpha(e)\leq \delta)=\delta\alpha^{-2}$. We
apply this estimate
with $\delta=1/k$ ($k$ large enough) and $r=j/k$ for all
$j=1,2,\ldots,Ak$. It follows that
$$E\Big[{\bf 1}_{E_A}\sum_{j=1}^\infty N^{1/k}_{j/k}(\alpha,c)\Big]\leq
2c\,A^2\,\alpha^{-2}.$$
In particular, if $E_k(\alpha,c,A)$ denotes the event
$\{\exists j\geq 1:N^{1/k}_{j/k}(\alpha,c)\geq 1\}\cap E_A$,
we have
$$P\Big[\liminf_{k\to \infty}E_k(\alpha,c,A)\Big]\leq
2c\,A^2\,\alpha^{-2}.\eqno{(6.1)}$$
Suppose that property (P) fails. Then, we may find
$t>\alpha$ and $s_1<s_2$ such that
$\wt \beta_{s_1}=\wt \beta_{s_2}=t$ and $0<m(s_1,s_2)<t-\alpha$, and
furthermore
$x_{m(s_1,s_2)}(\wt W_{s_1}(m(s_1,s_2)))=0$. We take $j$
such that $j/k<m(s_1,s_2)\leq(j+1)/k$, and observe that $x_{j/k}(\wt
W_{s_1}(j/k))<c$
for all $k$ sufficiently large, by the joint continuity of densities. Hence
by considering the excursion of $\wt \beta$ above level $j/k$ that contains
$s_1$,
we see that $N^{1/k}_{j/k}(\alpha,c)\geq 1$ for all $k$ large. Therefore, if
$F(\alpha,A)$ denotes the event on which (P) fails, we have
$$P[F(\alpha,A)\cap E_A]\leq P\Big[\liminf_{k\to \infty}E_k(\alpha,c,A))
\Big]\leq 2c\,A^2\,\alpha^{-2}.$$
Since $c$ was arbitrary, we have $P[F(\alpha,A)\cap E_A]=0$,
which completes the proof. \hfill$\square$
\medskip
\noindent{\bf Proof of Theorem 1.3}. The representation formula for $\wt Y_t$
implies that
$$\supp\wt Y_t=\{\wt W_s:\wt \beta_s=t\}.$$
(Note that the set on the right hand side is closed, by the
continuity properties of $\wt W$.) Hence if $w_1$ and $w_2$
belong to $\supp \wt Y_t$ and $w_1\not =w_2$, we can find
$s_1$ and $s_2$ such that $\wt \beta_{s_1}=\wt \beta_{s_2}=t$,
and $\wt W_{s_1}=w_1$, $\wt W_{s_2}=w_2$. With no loss of
generality, we can assume $s_1<s_2$. We claim that
$$m(s_1,s_2)=\inf\{r\in[0,t]:w_1(r)\not =w_2(r)\}.\eqno{(6.2)}$$
The inequality $m(s_1,s_2)\leq\inf\{r\in[0,t]:w_1(r)\not =w_2(r)\}$ is
immediate from the snake property (when $m(s_1,s_2)=0$ there is nothing
to prove). On the other hand, if we assume that there is
a rational $r\in(m(s_1,s_2),t)$ such that $\wt W_{s_1}(r)=\wt W_{s_2}(r)$,
then the monotonicity property implies $\wt W_s(r)= \wt W_{s_1}(r)$
for every $s\in[s_1,s_2]$ such that $\wt \beta_s\geq r$. Hence,
$$X_r=\int_0^{\tilde \tau} d\wt L^r_s\,\delta_{\wt W_s(r)}
\geq \int_{s_1}^{s_2} d\wt L^r_s\,\delta_{\wt W_s(r)}
=(\wt L^r_{s_2}-\wt L^r_{s_1})\,\delta_{W_{s_1}(r)},$$
which gives a contradiction since
$\wt L^r_{s_2}-\wt L^r_{s_1}>0$ by standard properties of linear Brownian
motion.

\smallskip
{}From now on, we assume $m(s_1,s_2)>0$.
Note that we have also $m(s_1,s_2)<t$ since we assumed that $w_1\not =w_2$.
By Lemma 6.1, we have $x_{m(s_1,s_2)}(\wt W_{s_1}(m(s_1,s_2)))>0$.
By monotonicity (and the fact that the measure $X_r$ gives no mass
to singletons), we get for every $r\in(m(s_1,s_2),t)$,
$$\int_{s_1}^{s_2} d\wt L^r_s=X_r((\wt W_{s_1}(r),\wt W_{s_2}(r)))
=\int_{\wt W_{s_1}(r)}^{\wt W_{s_2}(r)} dz\,x_r(z),$$
and by the continuity of densities, it follows that
$$\lim_{r\downarrow m(s_1,s_2)}
{\wt W_{s_2}(r)-\wt W_{s_1}(r)\over \wt L^r_{s_2}-\wt L^r_{s_1}}=
x_{m(s_1,s_2)}(\wt W_{s_1}(m(s_1,s_2))) >0.\eqno{(6.3)}$$
Thanks to (6.3), the behavior of $w_1(r)-w_2(r)$
as $r\downarrow \gamma_{w_1,w_2}=m(s_1,s_2)$ is reduced to that
of $\wt L^r_{s_2}-\wt L^r_{s_1}$.
Write $s_0$ for the (unique) time in $(s_1,s_2)$
such that $\wt \beta_{s_0}=m(s_1,s_2)$.
Standard results on Brownian path decompositions
show that, for events that depend only on the
asymptotic $\sigma$-field at time $0$,
the processes
$\{\wt \beta_{s_0-u} - \wt \beta_{s_0}, u\in[0,s_0-s_1]\}$
and
$\{\wt \beta_{s_0+u} - \wt \beta_{s_0}, u\in[0,s_2-s_0]\}$
behave
as two independent 3-dimensional Bessel processes.
It follows from this and the Ray-Knight theorem that the process
$\delta \to \wt L_{s_2}^{m(s_1,s_2) + \delta} - \wt L_{s_1}^{m(s_1,s_2)+
\delta}$
has the same local path properties (for $\delta$ close to 0)
as the sum of two
independent squares of 2-dimensional Bessel processes,
which is the square
of a 4-dimensional Bessel process.
If $\delta \to R_\delta$
is the square of a 4-dimensional Bessel process,
the law of the iterated logarithm shows that
$$\limsup_{\delta \downarrow 0}
{R_\delta \over
2  \delta \log |\log \delta | } = 1.$$
On the other hand, from the well-known rate of escape
for Brownian motion in space (Theorem 6 in [DE] combined
with time-inversion), we have for $\alpha>0$,
$$\lim_{\delta \downarrow 0}
{R_\delta \over
\delta |\log \delta|^ {-1 -\alpha} } = \infty.$$
We have just argued that the same properties hold
if we replace $R_\delta$ with
$\wt L_{s_2}^{m(s_1,s_2) + \delta} - \wt L_{s_1}^{m(s_1,s_2)+ \delta}$.
This and (6.3) imply Theorem 1.3.
\hfill$\square$
\bigskip
\noi{\bf Appendix}
\bigskip
\noi{\bf Proof of Lemma 5.3}. For a fixed value of $z$, the
estimate of Lemma 5.3 follows from [P2]. As we need uniformity in $z$,
we will provide a detailed argument.
Recall the notation from Subsection 5.2, and especially
the conventions concerning the constant $\alpha$.
Recall that $\cal G$
denotes the graph of
$X$ and for every integer $A\geq 1$ consider the event
$$E_A=\{{\cal G}\subset [0,A]\times [-A,A]; \sup_{t\geq \alpha,y\in \R}
x_t(y)\leq A\}.$$
Note that $P[E_A]\uparrow 1$ as $A\uparrow\infty$. (We use assumption (H)
when $\alpha=0$.)
The key step of the proof is to show the following
inequality for all $t\geq \alpha$
and $z\in \R$,
$$P[\{|X_{t+\delta}((-\infty,z])-X_t((-\infty,z])|\geq \delta^{{1\over
2}-\eta}\}
\cap E_A]\leq C\,\exp(-\delta^{-\kappa})\eqno{(A1)},$$
where the constants $C$ and $\kappa>0$ may depend on $A$ but not on $t,z$
and $\delta$.
To prove (A1), we may apply the Markov property at time $t$ and
reduce the problem to the case $t=0$. More precisely it is enough to
consider a super-Brownian motion
$\Gamma=(\Gamma_t,t\geq 0)$ with initial value $\Gamma_0(dz)=g(z)dz$,
with a function $g$ bounded above by $A$ and such that $\int g(z)dz\leq
2A^2$, and to prove
that for every $\delta\in(0,1)$,
$$P[|\Gamma_{\delta}((-\infty,0])-\Gamma_0((-\infty,0])|\geq
\delta^{{1\over 2}-\eta}]\leq
C\,\exp(-\delta^{-\kappa}).\eqno{(A2)}$$
Let us first bound
$P[\Gamma_{\delta}((-\infty,0])\leq\Gamma_0((-\infty,0])-\delta^{{1\over
2}-\eta}]$.
We know that for every $\lambda>0$,
$$E[\exp(-\lambda\Gamma_\delta((-\infty,0]))]=\exp(-\langle
\Gamma_0,u_\delta\rangle),$$
where $u_t(z)$ solves the integral equation
$$u_t(z)+{1\over 2}E_z\Big[\int_0^t
u_{t-r}(B_r)^2\,dr\Big]=\lambda\,P_z[B_t\leq 0],$$
if $B$ is a linear Brownian motion started at $z$ under $P_z$. The integral
equation gives the bound
$$u_t(z)\geq \lambda\,P_z[B_t\leq 0]-{\lambda^2\over 2} t.$$
We use this bound in the following estimates,
$$\eqalign{
&P[\Gamma_{\delta}((-\infty,0])\leq\Gamma_0((-\infty,0])-\delta^{{1\over
2}-\eta}]\cr
&\leq \exp(-\lambda\,\delta^{{1\over 2}-\eta}+\lambda\,\Gamma_0((-\infty,0]))
\,E[\exp(-\lambda\Gamma_\delta((-\infty,0]))]\cr
&\leq \exp(-\lambda\,\delta^{{1\over 2}-\eta}+{\lambda^2\over 2}\delta\langle
\Gamma_0,1\rangle)\,\exp\Big(\lambda\big(\Gamma_0((-\infty,0])-\int dz
g(z)P_z[B_\delta\leq 0]
\big)\Big).
}$$
Note that for every $\varepsilon>0$,
$$\Big|\int_{-\infty}^0 dzg(z)-\int dz g(z)P_z[B_\delta\leq 0]\Big|\leq
C_\varepsilon
\delta^{{1\over 2}-\varepsilon},$$
with a constant $C_\varepsilon$ depending only on $\varepsilon$ and $A$. By
choosing
$\lambda=\gamma^{-{1\over 2}+\varepsilon}$ with $0<\varepsilon<\eta$, we
arrive at
the desired estimate for
$P[\Gamma_{\delta}((-\infty,0])\leq\Gamma_0((-\infty,0])-\delta^{{1\over
2}-\eta}]$.
Slightly different arguments apply to
$P[\Gamma_{\delta}((-\infty,0])\geq\Gamma_0((-\infty,0])+\delta^{{1\over
2}-\eta}]$.
In fact, it is easier to observe that
$$\eqalign{&P[\Gamma_{\delta}((-\infty,0])\geq\Gamma_0((-\infty,0])
+\delta^{{1\over 2}-\eta}]
\cr
&\leq P[\langle \Gamma_\delta,1\rangle\geq \langle \Gamma_0,1\rangle
+{1\over 2}
\delta^{{1\over 2}-\eta}]
+P[\Gamma_{\delta}((0,\infty))\leq\Gamma_0((0,\infty))-{1\over
2}\delta^{{1\over 2}-\eta}].
}\eqno{(A3)}$$
We have just shown how to bound the second term on the right
hand side of (A3).
As for the
first term, we need simply recall that $\langle\Gamma_t,1\rangle$ is a
Feller diffusion
and use the fact that for $\lambda\in(0,{2\over \delta})$
$$E[\exp(\lambda\langle\Gamma_\delta,1\rangle)]
=\exp\Big({\lambda\langle \Gamma_0,1\rangle\over 1-{1\over 2}\lambda
\delta}\Big).
\eqno{(A4)}$$
This immediately leads to the estimate needed to complete the
proof of (A2) and (A1).

{}From (A1) and the Borel-Cantelli lemma, we get
that a.s. there is an integer $n_0(\omega)$ such that, for every
$n\geq n_0$, for every $t\geq 0$ of the form $t=j2^{-n}$ and every $z\in\R$
of the form $z=k2^{-n}$, we have
$$|X_{t+2^{-n}}((-\infty,z])-X_t((-\infty,z])|\leq (2^{-n})^{{1\over
2}-\eta}.$$
Note that for every fixed $z$, the process $t\to X_t((-\infty,z])$ has
continuous sample paths a.s. (see e.g. Corollary 6 in [P2]). The proof of
Lemma 5.3 is easily
completed thanks to this observation, the preceding bound and the usual
chaining argument.
\hfill$\square$

\bigskip
\noi{\bf Proof of Lemma 5.8}. This is very similar to the proof of (A1) above.
Note that the process $(X^*_{t+r},0\leq r\leq \delta)$
is a super-Brownian motion started at $X^*_t$, which is simply the
restriction of $X_t$ to $[z-\gamma,z+\gamma]$. Thanks to this
observation and the definition of $E(\delta)$, we see that it is enough to
prove the following statement. Let $\Gamma=(\Gamma_r,r\geq 0)$
be super-Brownian motion
with initial value $\Gamma_0(dz)=g(z)dz$. Assume that the function
$g$ vanishes outside $[-\gamma,\gamma]$ and that $c\leq g(z)\leq c^{-1}$
and $|g(z)-g(z')|\leq |z-z'|^{{1\over 2}-\rho}$ for all
$z,z'\in [-\gamma,\gamma]$. Then,
$$P[\Gamma_\delta((-\infty,0])\geq \Gamma_0((-\infty,0])+\delta^{{3\over
4}-\eta''}]
\leq \ov{C}\exp(-\delta^{-\ov{\kappa}}),\eqno{(A5)}$$
where the constants $\ov{C}$ and $\ov{\kappa}$ depend only on
$c,\eta',\eta''$ and $\rho$.

In a way similar to (A3) we first write
$$\eqalign{&P[\Gamma_{\delta}((-\infty,0])\geq\Gamma_0((-\infty,0])
+\delta^{{3\over 4}-\eta''}]
\cr
&\leq P[\langle \Gamma_\delta,1\rangle\geq \langle \Gamma_0,1\rangle
+{1\over 2}
\delta^{{3\over
4}-\eta''}]
+P[\Gamma_{\delta}((0,\infty))\leq\Gamma_0((0,\infty))-{1\over
2}\delta^{{3\over
4}-\eta''}].
}$$
Thanks to (A4), we see that, for $\lambda<2/\delta$,
$$\eqalign{P[\langle \Gamma_\delta,1\rangle\geq \langle \Gamma_0,1\rangle
+{1\over 2}
\delta^{{3\over
4}-\eta''}]
&\leq \exp(-\lambda(\langle \Gamma_0,1\rangle+{1\over 2}\delta^{{3\over
4}-\eta''}))
E[e^{\lambda\langle \Gamma_\delta,1\rangle}]\cr
&= \exp(-{\lambda\over 2}\delta^{{3\over 4}-\eta''})\,
\exp\Big({\lambda^2\langle \Gamma_0,1\rangle\delta/2\over 1-{1\over
2}\lambda\delta}\Big).
}$$
Since $\langle \Gamma_0,1\rangle\leq 2c^{-1}\gamma=2c^{-1}\delta^{{1\over
2}-\eta'}$, we get
a bound of the desired form by taking $\lambda=\delta^{-{3\over
4}+\varepsilon}$
with $\eta''>\varepsilon>\eta'$.

For the other term, we proceed as in the proof of Lemma 5.3:
$$P[\Gamma_{\delta}((0,\infty))\leq\Gamma_0((0,\infty))-{1\over
2}\delta^{{3\over
4}-\eta''}]\leq
\exp(\lambda(\Gamma_0((0,\infty))-{1\over 2}\delta^{{3\over
4}-\eta''}))E[e^{-\lambda\Gamma_\delta((0,\infty))}],$$
and
$E[e^{-\lambda\Gamma_\delta((0,\infty))}]=\exp(-\langle\Gamma_0,u_\delta\rangle)
$,
with $u_\delta(y)\geq \lambda P_y[B_\delta>0]-{1\over 2}\lambda^2\delta$.
It follows
that
$$\eqalign{&P[\Gamma_{\delta}((0,\infty))\leq\Gamma_0((0,\infty))-{1\over
2}\delta^{{3\over
4}-\eta''}]\cr
&
\leq \exp(-{1\over 2}\lambda\delta^{{3\over
4}-\eta''}+{\lambda^2\over 2}\delta\langle
\Gamma_0,1\rangle)\,\exp\Big(\lambda\big(\int_0^\infty dz g(z)-\int dz
g(z)P_z[B_\delta>0]
\big)\Big)\cr
&\leq  \exp(-{1\over 2}\lambda\delta^{{3\over
4}-\eta''}+c^{-1}\lambda^2\delta\gamma)\,
\exp(4\lambda\gamma^{{3\over 2}-\rho}),
}$$
where in the last line we used our assumption that
$|g(z)-g(0)|\leq |z|^{{1\over
2}-\rho}$ to
bound
$\int_0^\infty dz g(z)-\int dz g(z)P_z[B_\delta>0]$.
In view of the assumptions of Lemma 5.8, we can now choose
$\lambda=\delta^{{3\over
4}-\varepsilon}$, with $\eta''>\varepsilon>{3\over 2}\eta'+{\rho\over 2}$, and
we arrive at a bound of the desired form. This completes the proof.
\hfill$\square$

\vskip 1cm
\bigskip
\centerline{\bf REFERENCES}
\bigskip
\bigskip

\item{[AT]} R. Adler and R. Tribe,
Uniqueness for a historical SDE with a singular interaction.
J. Theoret. Probab. 11, 515--533 (1998)

\item{[BHM]} K. Burdzy, R. Ho\l yst and P. March,
A Fleming-Viot particle representation of Dirichlet
Laplacian (preprint)

\item{[Da]} D.A. Dawson,
Measure-valued Markov processes.
Ecole d'Et\'e de Probabilit\'es de Saint-Flour XXI---1991,
1--260, Lecture Notes in Math., 1541, Springer, Berlin, 1993.

\item{[DF]} D.A. Dawson, and K. Fleischmann, A continuous
super-Brownian motion in a super-Brownian medium.
J. Theoret. Probab. 10, 213--276 (1997)

\item{[DIP]} D.A. Dawson, I. Iscoe, and E.A. Perkins, Super-Brownian
motion: path properties and hitting probabilities.
Probab. Th. Rel. Fields 83, 135--205 (1989)

\item{[DP]} D.A. Dawson, and E.A. Perkins,
Historical processes. Mem. Amer. Math. Soc. 93, no. 454 (1991)

\item{[De]} J.F. Delmas, Super-mouvement brownien avec catalyse.
Stochastics Stochastics Rep. 58, 303--347 (1996)

\item{[DGL1]} D. D\"urr, S. Goldstein, and J.L. Lebowitz,
Asymptotics of particle trajectories in infinite one-dimensional
systems with collisions. Comm. Pure Appl. Math. 38, 573--597 (1985)

\item{[DGL2]} D. D\"urr, S. Goldstein, and J.L. Lebowitz,
Self-diffusion in a nonuniform one-dimen-{\allowbreak}sional system of point
particles with collisions. Probab. Th. Relat. Fields 75, 279--290 (1987)

\item{[DE]} A. Dvoretzky, P. Erd\" os, Some problems on
random walk in space. Proc. Second Berkeley Symp.
on Math. Stat. and Probab., pp. 353-367.
University of California Press, Berkeley 1951

\item{[Dy]} E.B. Dynkin, An Introduction to Branching Measure-Valued
Processes. CRM Monograph Series Vol. 6. American Mathematical Society,
Providence 1994

\item{[EK]} S.N. Ethier and T.G. Kurtz,
Markov processes: Characterization and convergence,
Wiley, New York, 1986.

\item{[EP]} S.N. Evans and E.A. Perkins, Measure-valued branching
diffusions with singular interactions. Canadian J. Math. 46, 120--168 (1994)

\item{[G]} R. Gisselquist, A continuum of collision process limit
theorems. Ann. Probab. 1, 231--239 (1973)

\item{[H]} T.E. Harris, Diffusion with ``collisions'' between
particles. J. Appl. Probab. 2, 323--338 (1965).

\item{[Ho]} R. Holley, The motion of a large particle.
Trans. Amer. Math. Soc. 144, 523--534 (1969)

\item{[KS]} N. Konno and T. Shiga,
Stochastic partial differential equations for some
measure-valued diffusions,
Probab. Th. Rel. Fields 79, 201--225 (1988)

\item{[L1]} J.F. Le Gall,
Marches al\'eatoires, mouvement brownien et processus de branchement.
S\'eminaire de probabilit\'es XXIII. Lecture Notes Math. 1372, pp. 258--274.
Springer 1989.

\item{[L2]} J.F. Le Gall,
Spatial Branching Processes, Random Snakes and Partial
Differential Equations, Lectures in Mathematics ETH Z\"urich,
Bikh\"auser 1999

\item{[NP]} J. Neveu, J. Pitman,
The branching process in a Brownian excursion.
S\'eminaire de probabilit\'es XXIII. Lecture Notes Math. 1372, pp. 248--257.
Springer 1989.

\item{[P1]} E. Perkins, Polar sets and multiple points for super-Brownian
motion. Ann. Probab. 18, 453--491 (1990)

\item{[P2]} E. Perkins,
On the continuity of measure-valued processes. In: Seminar
on Stochastic Processes 1990, pp. 261-268.
Progess in Probability 24. Birkh\"auser, Boston 1991

\item{[P3]} E. Perkins, On the martingale problem for
interactive measure-valued diffusions. Memoirs Amer. Math. Soc. 115, no. 549
(1995)

\item{[P4]} E. Perkins,
Dawson-Watanabe superprocesses and measure-valued diffusions.
Ecole d'Et\'e de Probabilit\'es de Saint-Flour XXIX---1999,
Lecture Notes in Math., Springer, to appear

\item{[R]} M. Reimers, One dimensional stochastic partial
differential equations and the branching measure diffusion.
Probab. Th. Rel. Fields 81, 319--340 (1989)

\item{[RY]} D. Revuz, M. Yor, Continuous Martingales and
Brownian motion. Springer, Berlin 1991

\item{[S]} F. Spitzer, Uniform motion with ellastic collision
of an infinite particle system. J. Math. Mech. 18, 973--989 (1968/69)

\vskip1truein

\parskip=0pt
\hbox{\vbox{
\obeylines
Krzysztof Burdzy
Department of Mathematics
University of Washington
Box 354350
Seattle, WA 98195-4350, USA
e-mail: burdzy@math.washington.edu
}
\hskip-8truecm
\vbox{
\obeylines
Jean-Fran\c cois Le Gall
DMA --- Ecole Normale Sup\'erieure
45, rue d'Ulm
75230 Paris Cedex 05, France
e-mail: legall@dma.ens.fr
{\ }
}
}

\bye